\documentclass[a4paper]{article}
\usepackage{epsfig}
\usepackage{amsmath}
\usepackage{amsfonts,amssymb}
\usepackage{amsthm}
\renewcommand{\ge}{\varepsilon}
\newcommand{\RR}{\mathbb{R}}

\newcommand{\CC}{\mathbb{C}}

\newcommand{\action}{\mathcal{A}}

\newcommand{\al}{\alpha}
\newcommand{\elli}{\mathcal{E}}

\let\R=\RR
\let\de=\partial

\DeclareMathOperator{\supp}{supp}
\newcommand{\bk}{\mathbf{k}}
\newcommand{\bn}{\mathbf{n}}
\newtheorem{theorem}{Theorem}[section]
\newtheorem{lemma}[theorem]{Lemma}

\newtheorem{definition}[theorem]{Definition}
\newtheorem{proposition}[theorem]{Proposition}
\newtheorem{remark}[theorem]{Remark}

\begin{document}


\title{Morse index properties of colliding solutions to the $N$-body
problem}

\author{Vivina Barutello\thanks{Partially supported by M.I.U.R., national project
    \textit{Metodi variazionali ed equazioni differenziali non lineari} and by I.N.D.A.M.,
    \textit{Istituto Nazionale di Alta Matematica}.} and Simone Secchi\thanks{Partially supported by M.I.U.R., national project
    \textit{Metodi variazionali ed equazioni differenziali non lineari}.}\\ \small
Dipartimento di Matematica e Applicazioni,
 \\ \small  Universit\`a di Milano--Bicocca via Cozzi 53, I-20125 Milano.\\
  \small Email: \texttt{vivina.barutello@unimib.it, simone.secchi@unimib.it}.}
\date{}
\maketitle

\begin{abstract}
  We study a singular Hamiltonian system with an $\al$-homogeneous
  potential that contains, as a particular case, the classical
  $N$--body problem. We introduce a variational Morse--like index for
  a class of collision solutions and, using the asymptotic estimates
  near collisions, we prove the non-minimality of some special classes
  of colliding trajectories under suitable spectral conditions
  provided $\al$ is sufficiently away from zero. We then prove some
  minimality results for small values of the parameter $\al$.

\end{abstract}

\section{Introduction}

In this paper we consider the second order Hamiltonian system
\begin{equation} \label{eq:dyn_syst}
 M \ddot x = \nabla U (x)
\end{equation}
with
\begin{equation}
\label{eq:pot_hom}
U(x) := \sum_{\substack{i,j=1 \\ i < j}}^N U_{ij}(x_i-x_j), \quad
U_{ij}(x_i-x_j) = \frac{m_i m_j}{|x_i-x_j|^\al},
\end{equation}
$x = (x_1,\dots, x_N) \colon [0,1] \to \R^{Nd}$, $d \geq 2$, $N \geq
2$, $\alpha \in (0,2)$, and $M = \operatorname{diag} [m_1,
\dots,m_N]$.  This system describes the well-known generalized
$N$--body problem, namely the motion of $N$ particles $x_1, \dots x_N$
of positive masses $m_1,\dots ,m_N$ under the external force $\nabla
U$ due to the generalized Kepler potential \eqref{eq:pot_hom}.  The
classical Keplerian case corresponds to the value $\al=1$. It is a
classical result, see \cite{painleve,vonzeipel}, that if $x$ is a
solutions of \eqref{eq:dyn_syst} on $[0,1)$ and if $x$ cannot be
extended to the whole interval $[0,1]$, then $\lim_{t \to 1}
U(x(t))=+\infty$; moreover, if $\|x\|$ remains bounded, then there
must be a collision at $t=1$, i.e. there exist two different indices
$i$, $j$ with $|x_i(t)-x_j(t)| \to 0$ as $t \to 1$.

It is evident that \eqref{eq:dyn_syst} has a rather delicate
variational structure, since the Euler--Lagrange action functional
\begin{equation}\label{eq:19}
\action_N (x) = \frac{1}{2} \int_0^1 | \dot x(t)|^2 \, dt+ \int_0^1 U (x(t))
\,dt
\end{equation}
may blow up along orbits $x(\cdot)$ that approach the collision set
\begin{equation} \label{collision_set}
\Delta = \{ x \in \RR^{Nd} : \exists i \neq j, \ x_i=x_j \}.
\end{equation}
Several recent papers are concerned with existence and qualitative
properties of \emph{collisionless solutions}, i.e. solutions such that
$x(t) \notin \Delta$ for all $t$. A first approach to avoid collision
solutions is the introduction of the so-called \emph{strong force}
assumption $\alpha \geq 2$ (see \cite{G}).  This constraint makes it
possible to prove that the Palais--Smale condition holds and to find
non-collision solutions by some standard tool of Critical Point
Theory.  Unfortunately, the Keplerian case does not satisfies such a
condition, for this reason much attention has been paid to the
complementary case $\alpha \in (0,2)$.  The bibliography about this
problem is huge, concerning the variational approach we cite, among
others,
\cite{AC-Z,bara91,cotizelati,C-Z,MR95k:58135,ft,serter1,serter2,tv}.

\medskip

In this paper we deal with some variational properties of solutions to
\eqref{eq:dyn_syst} possessing an isolated collision. Roughly
speaking, we will give an estimate of a generalized Morse index by
means of the asymptotic behavior of such a solution near the
collision. It is known that the action functional lacks regularity at
collision orbits, so that the usual Morse index cannot be defined.
This problem was overcome in \cite{MR95k:58135} by the technique of
approximate solutions. One of the results of that paper is an upper
bound on the number of total collisions for periodic solutions that
can be suitably approximated in the $H^1$--sense by solutions
corresponding to a regularized potential. The proof relies on the
construction of suitable variations introduced in \cite{serter1}.
Later, Riahi (see \cite{riahi}) generalized this result to solutions
with partial collisions, essentially by using the same method. We also
cite the paper \cite{tanaka}, where the author proves the existence of
one classical periodic solution in the case $1<\alpha <2$ and of one
generalized periodic solution with at most one collision in the case
$0<\alpha \leq 1$. The existence is proved again by the method of
approximate solutions, and the Author supplies some estimates on the
Morse index of these approximations. In the quoted papers, one of the
main ideas is that whenever the ratio between the dimension of the
space and the number of bodies involved in the collision is big
enough, then a collision gives a contribution to the Morse index of
the corresponding trajectory.

The main novelty of this paper consists in the use of the asymptotic
behavior near a collision in order to give an estimate on the Morse
index.  After fixing notation and reviewing some known facts (Section
2), we introduce in Section 3 the variational setting of our problem
and recall the main asymptotic estimates (see \cite{Sp,ft,Wi}) that
will be used to prove our main results.  In the next section we define
the generalized Morse index and provide in Theorem \ref{th:main} a
sufficient spectral condition on the asymptotic configuration,
ensuring that orbits with a single collision have an infinite index.
Theorem \ref{th:main} should be compared to the results of \cite{Del},
where non--minimality for a different class of colliding solutions is
proved. In particular, no condition like \eqref{eq:rel_eigen} appears
in that paper where, imposing reasonable assumptions about central
configurations and a symmetry assumption on the perturbation (the
Author perturbs the $N$-body potential with a term which is strongly
dominated in a $C^2$ sense near the collision set by the Newtonian
potential), it is shown that periodic orbits can be found with the
calculus of variations approach which avoid binary or triple
collisions. An additional assumption avoids total collapse orbits.

Section 5 is devoted to showing that the condition of the last section
is satisfied in some important cases, e.g. the collinear central
configuration of three equal masses or the polygonal configurations
for $N$ masses. In particular, we will show that our abstract theorem
applies for every $\alpha$ lying outside a small neighborhood of
$\alpha =0$.  Finally, Section 6 is somewhat complementary to the
previous ones.  Indeed, we analyze in deeper detail what happens in
the limit $\alpha \to 0$, and prove that under suitable assumptions,
families $(x_\alpha)_\alpha$ of one-collision solutions are
``minimal'', in the sense that the second derivative of the action
along compactly supported variations is positive.

\section{Preliminaries}
We consider the generalized keplerian potential defined in
\eqref{eq:pot_hom} on {\em collisionless configurations}
$x=(x_1,\ldots,x_N) \in \RR^{Nd} \setminus \Delta$, where $\Delta$ is
the collision set defined in \eqref{collision_set}.  The classical
keplerian interaction corresponds to the choice $\al = 1$.  We study
the dynamical system \eqref{eq:dyn_syst}, recalling that it is
conservative system, in the sense that the \textit{total energy}
\begin{equation} \label{eq:20} h = \frac12\langle
  M\dot{x},\dot{x}\rangle-U(x)= \frac12 \sum_{i=1}^N m_i|\dot x_i|^2 -
  U(x)
\end{equation}
is constant along solutions. Since the \emph{center of mass} moves
with a uniform motion, without loss of generality, we can fix it at
the origin, that is
\begin{equation*} 
\sum_{i=1}^N m_i x_i =0.
\end{equation*}
The potential \eqref{eq:pot_hom} will be then defined on the
configuration space
\begin{equation} \label{eq:21} \Lambda = \left\{x=(x_1,\ldots,x_N) \in
    \RR^{Nd} \setminus \Delta : \textstyle\sum_{i=1}^N
    m_ix_i=0\right\}
\end{equation}
For any $x \in \R^{Nd}$, the \textit{moment of inertia} is defined by
\begin{equation*}
I=I(x)=\langle Mx,x \rangle = \sum_{i=1}^N m_i |x_i|^2
\end{equation*}
and its gradient is simply
\begin{equation*}
\nabla I(x)=2Mx.
\end{equation*}
The following definition is quite standard.
\begin{definition} \label{def:cc} A {\em central configuration} is a
  critical point of the function $U$ constrained to the set $
  \mathcal{E} = \{ x \in \Lambda \mid I(x)=1 \}$. We will call
  $\mathcal{E}$ \textit{the standard ellipsoid}.
\end{definition}
\begin{remark} \label{rem:cc} An equivalent definition of central
  configuration is the following: $\bar x \in \mathcal{E}$ is a
  central configuration if there exists a solution of
  \eqref{eq:dyn_syst} of the form $x(t)=\phi (t) \bar x$, for some
  real-valued $C^2$-function $\phi$. For these classical facts, we
  refer to \cite{pollard,Wi}.
\end{remark}
Let us denote the radial and the angular components of $x \in \R^{Nd}
\setminus \{0\}$ by
\begin{equation*}
  r(x)=\sqrt{I(x)}, \qquad s(x) = \frac{x}{\sqrt{I(x)}}.
\end{equation*}
In particular $\mathcal{E}$ is now described by the simple condition $r=1$.
Since $U_{|\mathcal{E}}(s) = r^\al U(rs) = I^{\frac{\al}{2}}(x)U(x) $, it
follows easily that
\begin{equation} \label{eq:nabla_U} \nabla U_{|\mathcal{E}}(s)\cdot v =
  \frac{\al}{2}I^{\frac{\al}{2}-1}(x) U(x) \nabla I(x)\cdot v +
  I^{\frac{\al}{2}}(x) \nabla U(x)\cdot v,
\end{equation}
and
\begin{multline} \label{eq:nabla2_U} \nabla^2 U_{|\mathcal{E}}(s)(v,v)
  = I^{\frac{\al}{2}}(x) \nabla^2 U(x)(v,v) + \al
  I^{\frac{\al}{2}-1}(x) \langle\nabla I(x), v\rangle \langle\nabla
  U(x),v\rangle \\ {} + \frac{\al}{2}\left(\frac{\al}{2}-1 \right)
  I^{\frac{\al}{2}-2}(x)U(x)\langle\nabla I(x), v\rangle^2 +
  \frac{\al}{2} I^{\frac{\al}{2}-1}(x)U(x)\nabla^2 I(x)(v,v).
\end{multline}
As a consequence, when $x$ is a central configuration, that is when
$x\in \elli$ and $\nabla U_{|\mathcal{E}}(x)\cdot v = 0$ for every $v
\in T_x\elli$, we deduce from (\ref{eq:nabla_U}) that $\nabla U(x) =
-\al U(x)Mx$.  Replacing in (\ref{eq:nabla2_U}) we get
\begin{equation*}
\nabla^2 U_{|\mathcal{E}}(s)(v,v) = \nabla^2 U(x)(v,v) - \al(\al+2)
U(x)\langle Mx,v\rangle^2 + \al U(x)\langle Mv,v\rangle.
\end{equation*}
Since $v \in T_s \mathcal{E}$, we must have $\langle v , Mx \rangle
=0$, therefore the expression for the second derivative of
$U_{|\mathcal{E}}$ evaluated at $s$ is, for any $v \in T_{s}
\mathcal{E}$ with $\sum_i m_i v_i =0$,
\begin{equation}
\label{eq:nabla2_U_centr}
\nabla^2 U_{|\mathcal{E}}(s)(v,v) = \nabla^2U(x)(v,v) + \al U(x)\langle
Mv,v\rangle,
\end{equation}
where
\begin{equation}\label{eq:pot_hom_secondo}
\nabla^2 U(x)(v,v) = \al \sum_{i < j} m_i m_j \left[ (\al+2)
\frac{\langle x_i-x_j,v_i-v_j\rangle^2}{|x_i-x_j|^{\al+4}}
\right. \left. - \frac{|v_i-v_j|^2} {|x_i-x_j|^{\al+2}} \right]
\end{equation}
is the Hessian of $U$ on the whole space $\RR^{Nd}\setminus \Delta$.
When each $v_j$ is orthogonal to the vector space generated by $\{x_1,
\ldots ,x_N\}$, we deduce from~\eqref{eq:pot_hom_secondo} that the
Hessian of the potential $U$ is simply
\begin{equation}\label{eq:der_sec_normale}
\nabla^2 U(x)(v,v) = - \al \sum_{i < j} m_i m_j
\frac{|v_i-v_j|^2}{|x_i-x_j|^{\al+2}} = -\alpha \langle v,Av\rangle,
\end{equation}
where
\begin{equation}
  \label{eq:matrix_A} 
A(x)=[a_{ij}(x)], \quad a_{ij}(x)= \left \{
       \begin{array}{ll}
        \displaystyle \sum_{k \neq i} \frac{
        m_k}{|x_i-x_k|^{\al+2}} & i=j \\ \displaystyle -\frac{
        m_j}{|x_i-x_j|^{\al+2}} & i \neq j .
         \end{array}
     \right.
\end{equation}
\begin{remark}
Every tangent vector $v \in T_s \mathcal{E}$ can be seen as an $N$-uple of vectors $(v_1,\ldots,v_N)$, where each $v_j$ stands for the position of the $j$-th particle in the euclidean space $\R^d$. This justifies the slight abuse of looking at the Hessian $\nabla^2 U_{|\mathcal{E}}(s)$ as a quadratic form on $\R^N$. 
\end{remark}
In this case the expression for the constrained second derivative
\eqref{eq:nabla2_U_centr} can be recast as
\begin{equation} \label{eq:2nd_constr}
\nabla^2 U_{|\mathcal{E}}(s)(v,v) = \al \left( - \langle v,MAv \rangle +
  U(x)\langle Mv,v\rangle \right)
\end{equation}
for all $v=(v_1, \dots, v_N) \in \R^N$ with $\sum_{i=1}^N m_i v_i
=0$,

\section{The variational setting} \label{sec:vs}
It is well known that standard Critical Point Theory cannot be applied
to find solutions of \eqref{eq:dyn_syst} possessing a collision.
Indeed, the presence of collisions along a trajectory makes the action
function $\mathcal{A}_N$ (see definition \eqref{eq:19} below) possibly meaningless. As such, it might
even be impossible to say that a collision solution is a critical point of $\mathcal{A}_N$. 
For this reason, let us define the function spaces
\begin{equation*}
\Omega   = H^1\left((0,1),\Lambda\right), \qquad 
{\cal X} = H^1\left((0,1), \bar\Lambda \right),
\end{equation*}
where $\bar \Lambda = \{ x \in \R^{Nd} \mid \sum_{i=1}^n m_i x_i =0
\}$ is the closure of the set $\Lambda$ defined in \eqref{eq:21}.  The
elements of $\Omega$ will be termed \textit{collisionless orbits} and
their center of mass lies at the origin at every time.  Since each
element of $\Omega$ is a continuous function, it follows from standard
arguments that the action functional $\action_N \colon \Omega \to \R$
is smooth. Moreover, critical points of $\action_N$ inside the open
set $\Omega$ are collisionless, classical solutions of
\eqref{eq:dyn_syst}. It is clear that, in general, it is impossible to
extend the definition of $\action_N$ to ${\cal X}$, and it is
precisely this fact that prevents us from using standard tools for
studying colliding solutions to \eqref{eq:dyn_syst}.

In this paper we will take into account colliding solutions of
\eqref{eq:dyn_syst} with finite action and isolated collision. More
precisely, we introduce a class of ``good'' colliding solutions.
\begin{definition}\label{def:xcoll}
A {\em one-collision solution of \eqref{eq:dyn_syst}},
is a map $x\in C\left([0,1],\RR^{Nd}\right)\cap
C^2\left((0,1),\RR^{Nd}\right)$ such that
\begin{enumerate}
\item \label{item:2}
$x^{-1}(\Delta)\cap[0,1]=\{1\}$, i.e. no other collisions take
place in the time interval $[0,1)$;
\item $x$ solves pointwise the system \eqref{eq:dyn_syst} in the
interval $[0,1)$.
\item \label{item:1} $\action_N(x)<+\infty$;
\end{enumerate}
\end{definition}
\begin{remark}
Roughly speaking, a one-collision solution \textit{tends} to $\de\Lambda$ as,
and only as, $t \to 1$. Conditions \ref{item:2} and \ref{item:1} are strictly related 
(see for instance \cite{bftprogress,Sp}), we impose both not to enter in the details of this matter. 
\end{remark}
\begin{definition} \label{def:1cs}
Let $\bn := \{1,2,\ldots,N\}$. A \emph{colliding cluster for a
one-collision solution} $x$ is a subset $\bk \subset \bn$ such that
\begin{enumerate}
\item $x_i (1) = x_j (1)$ for all indices $i\neq j$ in $\bk$;
\item $x_i (1) \neq x_j (1)$ for all $i\in \bk$ and
      $j\in \bn \setminus \bk$.
\end{enumerate}
A collision will be termed \emph{total} if its associated cluster
$\bk=\bn$.
\end{definition}
The main property of a one-collision solution $x$ is that the action
$\action_N$ has directional derivatives at $x$ along
compactly--supported directions. This allows us to consider $x$ as a
``critical point'' of $\action_N$. The proof of the next lemma follows
trivially from Definition \ref{def:xcoll}.
\begin{lemma}
  Let $x$ be a one-collision solution of \eqref{eq:dyn_syst}.  Then
\begin{equation*}
\left.\frac{d}{d\varepsilon}\right|_{\varepsilon=0} {\cal A}_N
(x+\varepsilon \psi)=0, \quad \forall \psi \in C^\infty_0((0,1)).
\end{equation*}
\end{lemma}
Consider a one-collision solution $\bar x$ with a colliding cluster
$\bk \subset \bn$. Without loss of generality, we can assume $\bk =
\{1,2,\ldots,k\}$, so that the $N-k$ last components $(\bar{x}_{k+1},
\dots, \bar{x}_{N})$ of the one-collision solution $\bar x$ are kept fixed.
We define the restriction of the action functional
\begin{equation*} 
\action_{N,{\bf k}} (x) = \action_k(x) + \int_0^1 \Big[ \sum_{j=k+1}^N
|\dot{\bar x}_j|^2 + \frac{1}{2} \sum_{\substack{j_1,j_2 = k+1 \\ j_1
\neq j_2}}^N U_{j_1 j_2}(\bar x_{j_1} - \bar x_{j_2}) + \frac{1}{2}
\sum_{i=1}^{k}\sum_{j = k+1}^{N} U_{ji}(\bar x_{j} - x_{i}) \Big],
\end{equation*}
for every~$x \in H^1((0,1),\RR^{kd})$ and  
\begin{equation*}
\action (x) := \action_k(x) + \int_0^1 W(t,x)dt,
\end{equation*}
where
\begin{equation}
\label{eq:W}
W(t,x) := \frac{1}{2} \sum_{i=1}^{k}\sum_{j = k+1}^{N} U_{ji}(\bar
x_j(t) - x_i)
\end{equation}
is defined on $[0,1] \times \RR^{kd}$.
Since we have supposed that the action is finite at one-collision solutions the term 
$\int_0^1 \frac{1}{2} \sum U_{j_1 j_2}(\bar x_{j_1}-\bar x_{j_2})$ is finite (and constant), hence 
$\action_{N,{\bf k}}$ and $\action$ differ only by a constant. Of course, when all the bodies collide 
 the two functionals coincide on $H^1((0,1),\RR^{kd})$.
In the sequel we will deal with the functional $\action$.
\begin{remark}
  Since at $t=1$ the bodies in the cluster $\bk$ do not collide with
  those in $\bn \setminus \bk$ (and no other collision occurs in
  $[0,1)$), there exists an open set $\mathcal{U} \subset \RR^{kd}$
  such that $\left( \bar x_1([0,1]),\ldots,\bar x_k([0,1])\right)$ 
  $\subset {\cal U}$ and $W \in C^2 ([0,1] \times {\cal U})$.
\end{remark}

For simplicity, we will write $\bar x = (\bar x_1, \ldots , \bar x_k)$. Indeed, the terms involving the remaining components are of class $C^2$.
We define the radial and ``angular'' variables in the colliding cluster $\bk$
\begin{equation}
\label{rs}
r := |\bar x| = I^{\frac{1}{2}}(\bar x) \in \RR \quad s := \frac{\bar
x}{|\bar x|}.
\end{equation}
Since we are dealing with a total collision solution,
the following condition on the variable $r$ holds:
\begin{equation}
\label{collision}
\lim_{t \to 1}r(t)=0 \mbox{ and } r(t)\neq 0,
\quad \forall t \in [0,1).
\end{equation}
Condition~(\ref{collision}) means that the particles in $\bk$ collide
in their center of mass when $t=1$ and they do not have any other
collisions in the interval $[0,1)$.  Since
\begin{equation*}
\bar x = rs, \quad \dot{\bar x} = \dot{r}s + \dot{s}r, \quad
|\dot{\bar x}|^2 = \dot{r}^2+r^2 |\dot{s}|^2,
\end{equation*}
we can write the action functional at $x = \bar x$ in terms of the new
variables $(r,s)$ as
\begin{equation}
\label{action_rs}
\action(r,s) := \int_0^{1}
\frac{1}{2}r^{-(2+\al)}\left(r^{\frac{2+\al}{2}}\dot{r}\right)^2
+ r^{-\al} \left( \frac{1}{2}|r^{\frac{2+\al}{2}}\dot{s}|^2 +
U(s) + r^\al W(t,rs) \right)dt
\end{equation}
We consider the time scaling
\begin{equation}
\label{time_riscal}
dt=r^{\frac{2+\al}{2}}d\tau,
\end{equation}
and in the sequel we will note with a dot `` $\dot{}$ '' the
derivative with respect to the variable $t$ and with a prime `` $'$ ''
the one with respect to $\tau$.  Replacing~(\ref{time_riscal})
in~(\ref{action_rs}) we then obtain
\begin{equation}
\action (r,s) = \int_0^{\tau^*} \frac{1}{2} \left(
   r^{-\frac{2+\al}{4}} r'\right)^2 + r^{\frac{2-\al}{2}} \left(
   \frac{1}{2}|s'|^2 + U(s) + r^\al W\left(\int_0^{\tau}
   r^{\frac{2+\al}{2}},rs\right)\right) d\tau,
\label{action_rs_tau}
\end{equation}
where
\begin{equation*}
\tau^* = \int_0^{1} r^{-\frac{2+\al}{2}}dt.
\end{equation*}

In equation~(\ref{action_rs_tau}) we make the variable change
\begin{equation}
\label{rho}
\rho = r^{\frac{2-\al}{4}}, \quad \rho' = \frac{2-\al}{4}
r^{-\frac{2+\al}{4}} r'
\end{equation}
to obtain the action functional depending on $(\rho,s)$
\begin{equation}
\action (\rho,s) = \int_0^{\tau^*} \frac{1}{2}\left(
    \frac{4}{2-\al} \right)^2 (\rho')^2 + \rho^2 \left(
    \frac{1}{2}|s'|^2 + U(s) + \rho^{\frac{4\al}{2-\al}}
    W\left(\int_0^{\tau}
    \rho^\beta,\rho^{\frac{4}{2-\al}}s\right)\right) d\tau,
\label{action_rhos_tau}
\end{equation}
where 
\begin{equation}
  \label{eq:beta}
\beta := \frac{2(2+\al)}{2-\al} > 2.
\end{equation}
\begin{remark}
We notice that a set of quite similar variables were used in \cite{mcgehee} to study the dynamical system \eqref{eq:dyn_syst} from a geometrical viewpoint. As far as we know, the use of these coordinate in a variational setting is new.
\end{remark}
In the variables $(\rho,s,\tau)$, the Euler--Lagrange equations read
\begin{multline} \label{eq:2}
  -\left( \frac{4}{2-\al} \right)^2 \rho'' + \rho \left(|s'|^2 +
  2U_\al (s) \right) + \beta \rho^{\beta-1} \int_\tau^{\tau^*} \left[
    \rho^\beta \nabla_t W \left( \int_0^u
      \rho^\beta,\rho^{\frac{4}{2-\al}}s \right) \right] \, du \\
  {}+\beta \rho^{\beta-1} W \left( \int_0^\tau
    \rho^\beta,\rho^{\frac{4}{2-\al}}s \right) + \frac{4}{2-\al}
  \rho^{\frac32 \beta} \nabla W \left( \int_0^\tau
    \rho^\beta,\rho^{\frac{4}{2-\al}}s \right) \cdot s = \lambda_1
  \rho^{\beta-1},
\end{multline}
\begin{equation} \label{eq:4} 
  -2\rho \rho' s' - \rho \sp 2 s'' + \rho
  \sp 2 \nabla {U_\al}_{|\elli} (s) + \rho \sp {2\frac{4+\al}{2-\al}}
  \nabla_t W \left( \int_0^\tau \rho^\beta,\rho^{\frac{4}{2-\al}}s
  \right) = \lambda_2 s.
\end{equation}
It will be useful to have a more explicit formula for the Lagrange
multipliers $\lambda \sb 1$ and $\lambda_2$. We suppose here $W=0$,
since the computation is exactly the same in the general case. First of
all, we observe that the total (constant) energy $h$ (see
\eqref{eq:20}) can be written as
\begin{equation*}
h = \frac{1}{\rho^\beta} \left( \frac12 \left(
    \frac{4}{2-\al} \right)^2 (\rho')^2 + \rho^2 \left( \frac12
    |s'|^2 - U_\al (s) \right) \right).
\end{equation*}
Therefore, from \eqref{eq:2}
\begin{eqnarray*}
\frac{d}{d\tau} \left( \frac12 \left( \frac{4}{2-\al} \right)^2
  (\rho')^2 + \frac{\lambda_1}{\beta} \rho^\beta \right) &=&
\left( \left( \frac{4}{2-\al} \right)^2 \rho'' + \lambda_1 \rho^{\beta-1} \right) \rho' \\
&=&  \rho \left( |s'|^2 + 2 U_\alpha (s) \right) \rho'
=\frac{d}{d\tau} \left( \frac{\rho^2}{2} \right) \left( |s'|^2 + 2
  U_\al (s) \right),
\end{eqnarray*}
which implies that
\begin{eqnarray*}
\frac{d}{d\tau} \left( \rho^\beta \left( h_\al + \frac{\lambda_1}{\beta} 
\right) \right) &=&  \frac{d}{d\tau} \left( \frac12 \left( \frac{4}{2-\al} \right)^2 (\rho')^2 + \frac{\lambda_1}{\beta} \rho^\beta \right) + \frac{d}{d\tau} \left(  \frac{\rho^2}{2} \left( |s'|^2 -2U_\alpha (s) \right) \right) \\
&=& |s'|^2 \frac{d}{d\tau} \rho^2 + \frac{\rho^2}{2} \frac{d}{d\tau} \left( |s'|^2 -2U_\alpha (s) \right) \\
&=& 2 \rho \rho' |s'|^2 + \rho^2 s' \cdot s'' - \rho^2 \nabla {U_\al}_{|\elli}
(s) \cdot s' =0.
\end{eqnarray*}
Hence the constant $h_\al + \frac{\lambda_1}{\beta}$ must be zero, i.e.
\begin{equation} \label{eq:7}
\lambda_1 = - \beta h_\alpha.
\end{equation}
As for $\lambda_2$, we take the inner product of \eqref{eq:4} with $s$
and deduce immediately
\begin{equation} \label{eq:8}
\lambda_2 = \rho^2 |s'|^2,
\end{equation}
since $|s|^2=1$, $s \cdot s' =0$ and $s \cdot s'' = - |s'|^2$.
The next result describes the behavior of the new variables
$(\rho,s,\tau)$ and of the potential $U$ near the collision time. We
do not give a proof here, but we refer to \cite{PhD_vi} for a variational proof of these results, which were already proved in a different way (\cite{ft,Sp,Su}).
\begin{proposition}[Asymptotic estimates]
Let $x$ be a one-collision solution, $\bk$ its colliding cluster and
$\rho$, $s$, $\tau$ be defined by \eqref{rs}, \eqref{time_riscal}, \eqref{rho}
\label{prop:asymptotic}
The following properties hold true:
\begin{itemize}
\item[(a)] $\tau^* = +\infty$;
\item[(b)] There exists $b>0$ such that $\displaystyle\lim_{\tau \to
+\infty}\frac{\rho'}{\rho} = -\frac{2-\al}{4}\sqrt{2b}; $
\item[(c)] $\displaystyle \lim_{\tau \to +\infty} U(s(\tau))=b$;
\item[(d)] $\displaystyle \lim_{\tau \to {+\infty}} \operatorname{dist}
           \left({\cal C}^b,s(\tau)\right)= \lim_{\tau \to
           {+\infty}} \inf_{\bar s \in {\cal C}^b}|s(\tau) - \bar s
           |=0$, where
\begin{equation*}
\mathcal{C}^b := \left\{ s : U(s)=b, \nabla U_{|\mathcal{E}}(s)=0
\right\}.
\end{equation*}
is the set of central configurations (for the potential $U$) at level
$b$;
\item[(e)] $\displaystyle \int_0^\infty \frac{\rho'}{\rho}|s'|^2 < +\infty$.
\item[(f)] $\lim_{\tau \to +\infty} |s'(\tau)| =0$.
\end{itemize}
\end{proposition}

Point \emph{(d)} of Proposition~\ref{prop:asymptotic} does not mean
that the variable $s$ converges to an element of the set
$\mathcal{C}^b$; in this section we will deal with those collision
solutions that admit a limiting central configuration.  This fact is
expressed as follows.
\begin{definition} \label{def:as}
We say that a one-collision solution $x$ is \emph{asymptotic to a
central configuration} $s_0$ if $\lim\limits_{\tau \to +\infty} s(\tau)=s_0$.
More generally, we say that the one-collision solution $x$ is \emph{asymptotic
to the set of central configurations} $\mathcal{C}^b$ if {(d)} in Proposition
\ref{prop:asymptotic} is verified.
\end{definition}

\section{A class of colliding motions with non--trivial Morse index}
\label{sec:main}

We now introduce a notion of Morse index for one collision solutions with respect 
to the angular variable $s$.
The idea is to use the fact that in the new coordinates set ($\tau$, $\rho$  and $s$)
the collision take place at $+\infty$.
%
\begin{lemma} \label{lem:diff2}
  Let $x$ be a one-collision solution, and let $v \in T_s \mathcal{E}$
  be a compactly--supported function. Then
\begin{equation}\label{eq_diff2}
\frac{\de^2 \action}{\de s^2}(\rho,s)(v,v)=
 \int_0^{+\infty} \rho^2 \left[ |v'|^2 + \nabla^2 U_{|\mathcal{E}}(s) (v,v)
\right]
+  \rho^{2 \frac{6+\al}{2-\al}} \frac{\de^2}{\de x^2} W \left(\int_0^{\tau}
    \rho^\beta,\rho^{\frac{4}{2-\al}}s\right) (v,v)\, d\tau.
\end{equation}
\end{lemma}
\begin{proof}
The proof relies on very standard arguments. Take formula
\eqref{action_rhos_tau} and observe that the term
\begin{equation*}
\int_0^{+\infty} \frac{1}{2}\left( \frac{4}{2-\al} \right)^2
(\rho')^2 + \rho^2 \frac{|s'|^2}{2} \, d\tau
\end{equation*}
represents $\int_0^1 |\dot x|^2\, dt$, which is clearly smooth for $x
\in H^1((0,1),\R^{kd})$.  Therefore, we need to show that the
functional
\begin{equation*}
\Psi \colon (\rho,s) \in [0,+\infty) \times \mathcal{E} \mapsto
\int_0^{+\infty} \rho^2 \left( U(s) + \rho^{\frac{4\al}{2-\al}}
  W\left(\int_0^{\tau}
    \rho^\beta,\rho^{\frac{4}{2-\al}}s\right)\right) d\tau,
\end{equation*}
is twice differentiable in the variable $s$ along compactly supported directions.
Tt
follows at once from Proposition~\ref{prop:asymptotic} that $\rho$
decays exponentially fast and both $U(s(\tau))$ and $W(\int_0^\tau
\rho^\beta, \rho^{4/(2-\al)}s)$ remain bounded as $\tau \to
+\infty$. Hence we can apply the Dominated Convergence Theorem to show
that
\begin{eqnarray*}
\frac{\de^2 \Psi}{\de s^2}(\rho,s)(v,v) &=& \left. \frac{d^2}{d
\varepsilon^2} \right|_{\varepsilon =0} \Psi \left( \rho,
\frac{s+\varepsilon v}{\sqrt{I(s + \varepsilon v)}} \right) \\ &=&
\int_0^{+\infty} \rho^2 \nabla^2 U_{|\mathcal{E}}(s) (v,v) +
\rho^{2 \frac{6+\al}{2-\al}}\, \frac{\de^2}{\de x^2} W \left(\int_0^{\tau}
\rho^\beta,\rho^{\frac{4}{2-\al}}s\right) (v,v)\, d\tau.
\end{eqnarray*}
\end{proof}
\begin{definition} \label{def:morse}
  Let $x = \rho^{4/(2-\al)} s$ be a one-collision solution. We
  define the \emph{collision Morse index} $m_c=m_c (\action,\rho,s)$
  of $\action$ at $(\rho,s)$ as the supremum of all integers $m$ such
  that there exist $m$ linearly independent functions
  $\psi_1,\ldots,\psi_m \in C_0^\infty (\R^{+},T_s \mathcal{E})$ with
  the property that $\frac{\de^2 \action}{\de s^2}(\rho,s)$ is
  negative definite on $\operatorname{span} \left\{
  \psi_1,\ldots,\psi_m \right\}$. More precisely,
  $\frac{\de^2 \action}{\de s^2}(\rho,s)(v,v)< 0$ for all $v \in
  \operatorname{span} \left\{ \psi_1,\ldots,\psi_m \right\}$.
  Moreover, we will also say that $\action_N$ has collision Morse
  index $m_c$ at $x$.
\end{definition}

Our aim is to show that, under a suitable assumption on the
eigenvalues of the Hessian $\nabla^2 U_{|\mathcal{E}}(s_0)$, a
one-collision solution asymptotic to $s_0$ cannot be locally minimal
for the action functional \eqref{action_rhos_tau}. This is the content
of the next theorem.
%
\begin{theorem} \label{th:main}
  Let $x\in\mathcal{X}$ be a one-collision solution of
  \eqref{eq:dyn_syst} asymptotic to a central configuration $s_0$.
  Then the collision Morse index of $\mathcal{A}_n$ at $(\rho,s)$ is
  infinite, provided that the smallest eigenvalue $\mu_1$ of $\nabla^2
  U_{|\mathcal{E}}(s_0)$ satisfies
\begin{equation}
\label{eq:rel_eigen} 
\mu_1 < - \frac{(2-\al)^2}{8} U(s_0).
\end{equation}
\end{theorem}
\begin{proof}
  We introduce the variables $\rho$, $s$ and $\tau$ defined in \eqref{rho}, and according 
  with Definition \ref{def:morse} we will show
  that there exist infinitely many linearly independent functions
  $w_1, w_2, \dots$ such that $\frac{\de^2 \mathcal{A}}{\de
  s^2}(\rho,s)(w_i,w_i)<0$ for every index~$i$. For any smooth,
  compactly supported function $v$ such that $v(\tau) \in T_{s(\tau)}
  \mathcal{E}$ for all $\tau \geq 0$, we set $w=\rho v$. Then $w'=\rho' v + \rho v'$ and
\begin{equation} \label{eq:rhovi}
\rho^2 |v'|^2 = |w'|^2 + |\rho'|^2 |v|^2 - 2 \rho' w' v = |w'|^2 +
\left( \frac{\rho'}{\rho} \right)^2 |w|^2 - 2 \frac{\rho'}{\rho} w w'.
\end{equation}
In terms of $w$, equation \eqref{eq_diff2} becomes
\begin{multline}
\label{eq_diff2_bis} \frac{\de^2 \mathcal{A}}{\de s^2}(\rho,s)(v,v)=
\int_0^{+\infty} \left[ |w'|^2 + \left( \frac{\rho'}{\rho} \right)^2
  |w|^2 -2 \frac{\rho'}{\rho}w'w + \nabla^2 U_{|\mathcal{E}}(s) (w,w)
  \right. \\ + \left. \rho^{2\beta} \frac{\de^2}{\de
    x^2} W \left( \int_0^\tau \rho^\beta,\rho^{\frac{4}{2-\al}}s
  \right)(w,w) \right]\, d\tau.
\end{multline}
Setting
\begin{multline*}
Q(w)=
\int_0^{+\infty} \left[ |w'|^2 + \left( \frac{\rho'}{\rho} \right)^2
  |w|^2 -2 \frac{\rho'}{\rho}w'w + \nabla^2 U_{|\mathcal{E}}(s) (w,w)
  \right. \\ + \left. \rho^{2\beta} \frac{\de^2}{\de
    x^2} W \left( \int_0^\tau \rho^\beta,\rho^{\frac{4}{2-\al}}s
  \right)(w,w) \right]\, d\tau.
\end{multline*}
we will prove that $Q<0$ on a vector space of infinite dimension.
Let $0< \ell_1 < \ell_2$ be arbitrary numbers, and take a positive
real function $\varphi \in C_0^\infty (\ell_1,\ell_2)$; let
$\{\tau_n\}_n$ be a strictly increasing, divergent sequence of
positive numbers.  We define $w_n(\tau)= \varphi(\tau - \tau_n)\xi$,
where $\xi\in T_{s_0} \mathcal{E}$ will be chosen hereafter. In
particular $w_n \in C_0^\infty
((\ell_1+\tau_n,\ell_2+\tau_{n}),T_{s_0}\mathcal{E})$.  It follows
from Proposition \ref{prop:asymptotic} that the following estimates
hold:
\begin{equation*}
\int_0^{+\infty} \rho^{2\beta} \frac{\de^2}{\de x^2}
  W \left( \int_0^\tau \rho^\beta,\rho^{\frac{4}{2-\al}}s
  \right)(w_n,w_n)\, d\tau \leq C_1
  \int_{\ell_1+\tau_n}^{\ell_2+\tau_{n}}
  \rho^{2\beta} |w_n|^2 \, d\tau \leq C_1
  \|\varphi\|_\infty e^{-C_2 \tau_n},
\end{equation*}
\begin{equation*}
\int_{0}^{+\infty} \frac{\rho'}{\rho} w_n w_n'\, d\tau
=\int_{\ell_1+\tau_n}^{\ell_2+\tau_n} \frac{\rho'}{\rho} w_n w_n' \,
d\tau = -\frac{2-\al}{4} \sqrt{2 U(s_0)}
\int_{\ell_1+\tau_n}^{\ell_2+\tau_n} w_n w_n'\, d\tau + o(1) =o(1)
\end{equation*}
as $n \to +\infty$. In a similar way,
\begin{equation*}
\int_0^{+\infty} \left| \frac{\rho'}{\rho} \right|^2 |w_n|^2 \, d\tau
= \frac{(2-\al)^2}{8} U(s_0) \int_{\ell_1+\tau_n}^{\ell_2+\tau_n}
|w_n|^2\, d\tau + o(1)
\end{equation*}
Putting together these estimates, using the continuity of $\nabla^2 U$
at $s_0$ and the fact that $x$ is asymptotic to the central
configuration $s_0$, we get
\begin{equation}\label{eq:penultima}
Q(w_n)=
\int_{\ell_1+\tau_n}^{\ell_2+\tau_n} \left( |w_n'|^2 +
\frac{(2-\al)^2}{8} U(s_0) |w_n|^2 + \nabla^2
U_{|\mathcal{E}}(s_0)(w_n,w_n) \right) \, d\tau + o(1).
\end{equation}
We choose now $\xi \in T_{s_0} \mathcal{E}$ as a normalized eigenvalue of
$\nabla^2 U_{|\mathcal{E}}(s_0)$ corresponding to the eigenvalue $\mu_1$,
verifying assumption \eqref{eq:rel_eigen}. Then equation
\eqref{eq:penultima} becomes
\begin{equation*}
Q(w_n)=
\int_{\ell_1+\tau_n}^{\ell_2+\tau_n} \left[ |w_n'|^2 + \left(
\frac{(2-\al)^2}{8} U(s_0) + \mu \right) |w_n|^2 \right] \, d\tau + o(1).
\end{equation*}
Since \emph{both} $\varphi$ and its support $(\ell_1,\ell_2)$ are
arbitrary, it follows from \eqref{eq:rel_eigen} that $\frac{\de^2
\mathcal{A}}{\de s^2}(\rho,s)(w_n,w_n)<0$ for all $n \gg 1$. We can
now repeat the same construction with different choices of $\varphi$
and of the sequence $\{\tau_n\}_n$, and build a countable family of
functions $\{w_n\}_n$ with disjoint supports and such that
$\frac{\de^2 \mathcal{A}}{\de s^2}(\rho,s)(w,w)<0$ for all $w \in
\operatorname{span} \{w_1,w_2,\dots\}$. From Definition
\ref{def:morse} it follows that the collision Morse index of $\action$
at $x$ is infinite.
\end{proof}
In the next section we will present some concrete examples in which
our Theorem~\ref{th:main} applies.

\section{Applications of Theorem \ref{th:main}}
\label{sec:appl}
In this section we discuss the applicability of Theorem~\ref{th:main} to concrete examples of limiting central
configurations. Clearly, the hardest assumption to check is
inequality~\eqref{eq:rel_eigen}. Since it is known that the smallest
eigenvalue $\mu_1$ of $\nabla^2 U_{| \mathcal{E}} (s_0)$ at the
central configuration $s_0$ is characterized by
\begin{equation}
\mu_1 = \min \left\{ \nabla^2 U_{|\mathcal{E}}(s_0)(v, v) \mid v \in
T_{s_0} \mathcal{E}, \; \sum_i m_i v_i=0, \; \|v \| = 1 \right\} ,
\end{equation}
using \eqref{eq:nabla2_U} we obtain that \eqref{eq:rel_eigen} is implied by the existence of a vector
$v \in T_{s_0} \mathcal{E}$ such that $\sum_i m_i v_i=0,
\; \|v\| = 1$ and
\begin{equation} \label{eq:33}
 \nabla^2 U(s_0)(v,v) + \alpha U(s_0) \langle Mv,v \rangle < - \frac{(2-\alpha)^2}{8} U(s_0).
\end{equation}
In particular, when all the masses are equal to 1, we obtain the simpler condition
\begin{equation} \label{eq:nuova}
\nabla^2U_{|{\cal E}}(s_0)(v,v) < -\frac{(2-\al)^2}{8} \,U(s_0).
\end{equation}
When each $v_j$ is assumed to be orthogonal to the vector space
generated by the configuration $s_0$ (see~\eqref{eq:der_sec_normale})
we can introduce the square matrix $A$, defined
in~\eqref{eq:matrix_A}. Hence \eqref{eq:33} and \eqref{eq:nuova}
are satisfied provided we can find a vector $v =
(v_1,\dots,v_N) \in \R^N$, such that $\|v\|=1$ and
\begin{equation}
\label{eq:2555}
-\langle v,MAv \rangle + U(s_0)\langle v,Mv \rangle
< -\frac{(2-\al)^2}{8\al} \,U(s_0)
\end{equation}
and
\begin{equation}
\label{eq:25} 
\langle w,Aw \rangle > \frac{(\al+2)^2}{8\al} \,U(s_0),
\end{equation}
respectively. We will prove that for a wide range of values of~$\al$
(including the value $\al=1$) the collinear central
configuration of three equal masses and the regular $N$-gon
configuration satisfy \eqref{eq:rel_eigen}, showing \eqref{eq:25}.
In particular, in the second case, when $N$ is even, we will prove that 
\eqref{eq:25} is satisfied for a vector $w\in\RR^N$ that verifies the hip-hop symmetry 
(see \cite{cv} and \cite{tv}).

\subsection{Collinear central configurations for three equal masses}
We consider the collinear central configuration of three particles of
masses $m_1=m_2=m_3=1$, lying on a straight line
\begin{equation*}
s_0 = \left(( -1/\sqrt 2,0),(0,0),(1/\sqrt 2,0)\right).
\end{equation*}
We perform a planar variation as follows:
\begin{equation*}
v=\left((\cos \theta, \sin \theta),(0,-2 \sin \theta),(-\cos \theta,
\sin \theta)\right),
\end{equation*}
where $\theta \in [0,\pi /2]$. We remark that $v=(v_1,v_2,v_3) \in
T_{s_0}\mathcal{E}$, $\sum_{i=1}^3 v_i =0$, and $\|v\|^2 = 2
(1+2 \sin^2 \theta)$.  With these choices the Hessian at the
configuration $s_0$ is (see~\eqref{eq:pot_hom_secondo})
\begin{equation*}
\nabla^2 U(s_0)(v,v) = 2\al \left\{\cos^2\theta\left[
  2(\al+10)2^{\al/2} + (\al+1)2^{-\al/2}\right] -18\cdot
  2^{\al/2}\right\}
\end{equation*}
Therefore, after dividing out by $\|v\|^2$, \eqref{eq:rel_eigen}
reads
\begin{multline}
\label{eq:storto}
\frac{1}{1+2\sin^2\theta}\left\{\cos^2\theta\left[
  2(\al+10)2^{\al/2} + (\al+1)2^{-\al/2}\right] -18\cdot
  2^{\al/2}\right\} \\ {} < -\frac{(\al+2)^2}{8\al}
  \left(2\cdot 2^{\al/2} + 2^{-\al/2}\right)
\end{multline}
It is apparent that the most convenient choice, in order to get the
widest range of $\al$'s, is $\theta = \pi /2$, i.e. to take normal
variations. Hence \eqref{eq:storto} reduces to
\begin{equation}
\label{eq:dis_coll}
\frac{6 \cdot 2^{\al}}{2\cdot 2^\al + 1} >
\frac{(\al+2)^2}{8\al}
\end{equation}
Let $f(\al):=\frac{6 \cdot 2^{\al}}{2\cdot 2^\al + 1}$ and
$g(\al):=\frac{(\al+2)^2}{8\al}$ be respectively the left and
right hand side of \eqref{eq:dis_coll}; since $g$ strictly decreases
on $(0,2]$, $f$ strictly increases $[0,2]$ and $f(0)=g(6-4\sqrt{2})$,
we conclude the existence of $\bar \al < 6-4\sqrt{2}$ such that for
every $\al \in [\bar \al,2]$ the inequality \eqref{eq:dis_coll}
holds true.

In a similar way, we can consider the central configuration of three
masses $m_1=m_3$ but $m_2$ different. Indeed, in this case we have a
central configuration $s_0$ whose points are
\begin{equation*}
s_0=\left((-1/\sqrt{2m_1},0),(0,0),(1/\sqrt{2m_1},0)\right).
\end{equation*}
We then choose again the normal variation
\begin{equation*}
v_1=(0,1), \quad v_2=(0,-2), \quad v_3=(0,1)
\end{equation*}
and observe that condition \eqref{eq:2555} reads now
\begin{equation*}
\frac{9\cdot 2^{\frac{\al+2}{2}} m_1^{\frac{\al+4}{2}} m_2 -
(m_1+2m_2) \left( 2^{\frac{\al+2}{2}} m_1^{\frac{\al+2}{2}} m_2
+ 2^{-\frac{\al}{2}} m_1^{\frac{\al+4}{2}}
\right)}{2^{\frac{\al+2}{2}} m_1^{\frac{\al+2}{2}} m_2 +
2^{-\frac{\al}{2}} m_1^{\frac{\al+4}{2}}} >
\frac{3(2-\al)^2}{8\al}
\end{equation*}
After some simplifications, this is equivalent to the inequality
\begin{equation*}
\frac{9 \cdot 2^{\al+1} m_1 m_2 - (m_1 +2m_2)(2^{\al+1}
m_2+m_1)}{2^{\al+1} m_2 + m_1} > \frac{3(2-\al)^2}{8\al}.
\end{equation*}
Since this inequality is homogeneous with respect to the masses $m_1$
and $m_2$, we can suppose now $m_2=1$. Hence we should solve
\begin{equation} \label{eq:dueuguali}
\frac{2^\al (16m_1+4) -m_1^2
+2m_1}{2^{\al+1}+m_1}>\frac{3(2-\al)^2}{8\al}.
\end{equation}
Set now
\begin{equation*}
f(\al) = \frac{2^\al (16m_1+4) -m_1^2 +2m_1}{2^{\al+1}+m_1},
\quad g(\al) = \frac{3(2-\al)^2}{8\al}.
\end{equation*}
One checks easily that $g(2)=0$ and $g$ is a positive, strictly
decreasing function on $(0,2)$. Moreover, since
\begin{equation*}
\frac{f'(\al)}{2^\al \log 2} = \frac{18m_1^2}{(2^{\al+1}+m_1)^2},
\end{equation*}
the function~$f$ is strictly increasing to the value $f(2) =
\frac{-m_1^2+66m_1+16}{m_1+8}$. We conclude that we can find a number
$\al^*>0$ such that \eqref{eq:dueuguali} is satisfied for all
$\al > \al^*$ if and only if $f(2)>g(2)$, i.e. $m_1^2-66m_1-16<0$, or 
$m_1 < 33+\sqrt{33^2+16}$. The
Newtonian case $\al=1$ is admissible if and only if $f(1)>g(1)$,
i.e.  $\frac{2(16m_1+4)-m_1^2+2m_1}{4+m_1} > 3/8$, or $8m_1^2-269m_1 + 52 >0$. Hence $m_1 < 34$ is enough.


\begin{remark}
For the collinear configuration of three equal masses we can try to
verify \eqref{eq:rel_eigen} instead of the stronger \eqref{eq:nuova}.
Observe that the vector $(1,1,1)$ is an eigenvector for $A$ with
eigenvalue $0$. Hence we restrict the matrix $A$ to the space
orthogonal to this vector that is $Y=\{ v=(v_1,v_2,v_3):\sum_i
v_i=0 \}$ spanned by
\begin{equation*}
  w_1=(1,0,-1), \quad w_{2}=(0,1,-1).
\end{equation*} 
If $B=[b_{hk}]$ denotes the symmetric matrix $A$ restricted to the
space $Y$ we have that
\begin{equation*}
b_{hk} = w_h^T A w_k = a_{hk} - (a_{h3}+a_{k3})+a_{33} = b_{kh},
\end{equation*}
where $w^T$ denotes transposition of the vector $w$. Explicitly,
\begin{equation*}
B= \begin{pmatrix}
2\gamma + 4\gamma^{-1} & \gamma+2\gamma^{-1} \\ \gamma+2\gamma^{-1} &
    5\gamma + \gamma^{-1}
\end{pmatrix},
\quad \hbox{with $\gamma = 2^\frac{\al+2}{2}$},
\end{equation*}
and its eigenvalues are
\begin{equation*}
  \lambda_B^{1,2} = \frac{7\gamma +5\gamma^{-1} \pm \sqrt{13\gamma^2
  -2 +25\gamma^{-2}}}{2}.
\end{equation*}
It is easy to check that condition \eqref{eq:25} is implied by the
inequality
\begin{equation}
\label{eq:condition}
\frac{7\gamma +5\gamma^{-1} + \sqrt{13\gamma^2 -2 + 25\gamma^{-2}}}{2}
> \frac{(2+\al)^2}{8\al} \left( \gamma + 2\gamma^{-1} \right).
\end{equation}
This approach gives a wider range of ``good'' values for the parameter
$\al$ (as we can see in Figure \ref{fig:trecorpi}), but is clearly
impossible to apply in more general situations.
\begin{figure}[th!]
\begin{center}
  {\psfig{figure=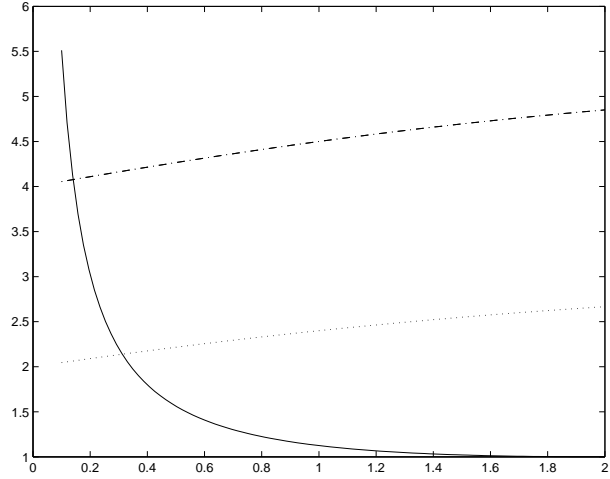,width=8.0cm}}
\end{center}
\caption{\label{fig:trecorpi} Validity of \eqref{eq:dis_coll} in the
          dotted line compared with the validity of
          \eqref{eq:condition}. The continue line represents the
          quantity $\frac{(2+\al)^2}{8\al}$.}
\end{figure}
\end{remark}

\subsection{The regular $N$-gon central configuration}

We now consider the case of a planar central configuration of $N \geq
4$ equal bodies with equal masses, lying at the vertices of a regular
$N$-gon inscribed in a circle of radius $1/\sqrt{N}$.  In the sequel,
we will systematically use the notation
\begin{equation*}
r_{ij} = r_{ij}(N) = \frac{2}{\sqrt{N}} \sin \left( \frac{|i-j|}{N}
\pi \right)
\end{equation*}
for the distance between the $i$-th and the $j$-th bodies. We remark
that $r_{ij}=r_{1k}$, where $k=|i-j|+1$.

The aim of this section is to show that when $s_0$ is the polygonal
central configuration, then relation \eqref{eq:rel_eigen} is verified
for a whole interval of $\al$'s including the Newtonian case
$\al=1$. This will prove that collisions ending up in a polygonal
configuration cannot be minima for the action.

Condition \eqref{eq:25} in this case reads
\begin{equation} \label{eq:25bis}
\sum_{i,j} a_{ij} w_i w_j > \frac{(\al+2)^2}{8\al} U(s_0) =
\frac{(\al+2)^2}{8\al} \frac{N}{2} \sum_{k=2}^N
\frac{1}{r_{1k}^\al}.
\end{equation}
We rewrite \eqref{eq:25bis} in the form
\begin{equation} \label{eq:36}
\Psi_N(\al) := \frac{2}{N} \frac{\sum_{i,j} a_{ij} w_i
w_j}{\sum_{k=2}^N r_{1k}^{-\al}} > \frac{(\al+2)^2}{8\al}.
\end{equation}
Define, for simplicity, $\tilde{r}_{ij} = \tfrac{\sqrt N}{2} r_{ij} =
\sin \left( \tfrac{|j-i|}{N}\pi \right)$, $\tilde{a}_{ij} = \left(
\tfrac{\sqrt 2}{N} \right)^{\al+2} a_{ij}$, so that the matrix $\tilde{A} = [\tilde{a}_{ij}]$ can be constructed by writing $\tilde{r}_{ij}$ instead of $r_{ij}$ in \eqref{eq:matrix_A}.
Observe that
\begin{equation*} 
\Psi_N(\al)= \frac12 \frac{\sum_{i,j} \tilde{a}_{ij} w_i
w_j}{\sum_{k=2}^N \tilde{r}_{1k}^{-\al}} = \frac12
\frac{\sum_{i=1}^N \tilde{a}_{ii} w_i^2 + 2 \sum_{i<j} \tilde{a}_{ij}
w_i w_j}{\sum_{k=2}^N \tilde{r}_{1k}^{-\al}}.
\end{equation*}
We choose $w$ as follows:
\begin{equation} \label{eq:w}
w=
\begin{cases}
(1/2,-1/2,1/2,-1/2) &\text{for $N=4$} \\
(1/\sqrt{2},-1/\sqrt{2},0,\dots,0) &\text{for $N \geq 5$}.
\end{cases}
\end{equation}
\begin{remark} \label{rem:dueindici}
Equivalently, when $N \geq 5$, we could choose $w$ such that $w_i =
1/\sqrt 2$, $w_{i+1}=-1/\sqrt 2$, and $w_k=0$ for any~$k \notin
\{i,i+1\}$. We will use this observation later on.
\end{remark}
Therefore
\begin{eqnarray*}
\Psi_4(\al) &=& \Phi_4 (\al) + \frac12 \frac{2
\tilde{r}_{12}^{-\al-2}-\tilde{r}_{13}^{-\al-2}}{\sum_{k=2}^4
\tilde{r}_{1k}^{-\al}}\\ \Psi_N(\al) &=& \Phi_N (\al) +
\frac12 \frac{\tilde{r}_{12}^{-\al-2}}{\sum_{k=2}^N
\tilde{r}_{1k}^{-\al}} \qquad (N \geq 5),
\end{eqnarray*}
where
\begin{equation}
\Phi_N (\al) := \frac12 \frac{\sum_{k=2}^N
\tilde{r}_{1k}^{-\al-2}}{\sum_{k=2}^N \tilde{r}_{1k}^{-\al}},
\quad\text{for all $N \geq 4$}.
\end{equation}
We now state and prove some technical lemmata that are useful for the
proof of the main theorem of this section. The first one is a simple
exercise in first year calculus.
\begin{lemma} \label{lem:scemo}
  Let $\{b_j\}_{j=1,\ldots,n}$ be a family of $n$ positive real
  numbers such that
\begin{equation*}
b_1 > b_2 > \ldots > b_n > 1.
\end{equation*}
Then the two functions
\begin{equation*}
f(x) := \frac{\sum_{j=1}^n b_j^{x+2}}{\sum_{j=1}^n b_j^{x}} \quad
\mbox{and} \quad g(x) := \frac{1+\sum_{j=1}^n
b_j^{x+2}}{1+\sum_{j=1}^n b_j^{x}}
\end{equation*}
are strictly increasing on the interval~$[0,2]$.
\end{lemma}

\begin{proof}
One computes the first derivatives of $f$ and $g$ by the rule for a
quotient. Then the sign of these derivatives is the sign of the
numerator.  When dealing with $f$, after simplifying some terms, we
end up with a sum of terms like
\begin{equation*}
b_i^x b_j^x (b_i^2 -b_j^2) (\log b_i - \log b_j), \quad (i<j)
\end{equation*}
and these are all positive, because of the monotonicity of the family
$\{b_j\}$. When dealing with $g$, an extra term appears due to the
presence of the number 1. Nevertheless, we easily check that the extra
term is just
\begin{equation*}
\sum_j b_j^x (b_j^2-1) \log b_j ,
\end{equation*}
which is positive since $b_j >1$.
\end{proof}

\begin{lemma} \label{lem:phi}
  For every $N \geq 4$, the function $\Phi_N \colon [0,2] \to \R$ is
  strictly increasing and satisfies
\begin{equation} \label{eq:28}
\Phi_N (0) > \frac{N-1}{N}.
\end{equation}
\end{lemma}
\begin{proof}
  The monotonicity of $\Phi_N$ follows easily from Lemma
  \ref{lem:scemo}, by exploiting the symmetry of the regular $N$-gon
  with respect to a straight line passing through a fixed vertex. If
  $N$ is odd, one has 
\begin{equation*}
\Phi_N (\al) := \frac12 \frac{\sum_{k=2}^{(N+1)/2}
\tilde{r}_{1k}^{-\al-2}}{\sum_{k=2}^{(N+1)/2} \tilde{r}_{1k}^{-\al}},
\end{equation*}
where $0<\tilde{r}_{1k} = \sin \left(\frac{k-1}{N}\pi \right) \leq \sin \left(\frac{N-1}{N}\pi \right) <1$.
Now we can use the monotonicity of $f$. If $N$ is even one
  uses the monotonicity of $g$ with $a_j :=\sin \left( \frac{j}{N} \pi
  \right)$, $j=1,\ldots,\left[ \frac{N+1}{2}\right]-1$.  The proof of
  \eqref{eq:28} is equivalent to the following inequality
\begin{equation*}
\sum_{k=2}^N \tilde{r}_{1k}^{-2} > \frac{2(N-1)^2}{N}
\end{equation*}
Since $t \mapsto t^{-2}$ is a convex function on $(0,\infty)$, the
discrete Jensen inequality tells us that
\begin{equation*}
 \bigg( \frac{1}{N-1} \sum_{k=2}^N \tilde{r}_{1k} \bigg)^{-2} \leq
\frac{1}{N-1} \sum_{k=2}^N \tilde{r}_{1k}^{-2}, \quad \hbox{i.e.}
\end{equation*}
\begin{equation*}
\frac{(N-1)^3}{\left( \sum_{k=2}^N \tilde{r}_{1k} \right)^2} \leq
\sum_{k=2}^N \tilde{r}_{1k}^{-2}.
\end{equation*}
Hence the following inequality implies \eqref{eq:28}
\begin{equation*}
\frac{(N-1)^3}{\left( \sum_{k=2}^N \tilde{r}_{1k} \right)^2} >
\frac{2(N-1)^2}{N},\quad \hbox{or}
\end{equation*}
\begin{equation*}
\sum_{k=2}^N \tilde{r}_{1k} < \sqrt{\frac{N(N-1)}{2}}.
\end{equation*}
Since $\tilde{r}_{1k} = \operatorname{Im} e^{i \frac{k-1}{N}\pi}$, we
can easily compute
\begin{equation*}
\sum_{k=2}^N \tilde{r}_{1k} = \operatorname{cotan} \frac{\pi}{2N} =
\frac{\sin \frac{\pi}{N}}{1-\cos \frac{\pi}{N}}
\end{equation*}
Therefore, we have to prove that
\begin{equation*} 
\frac{\sin \frac{\pi}{N}}{1-\cos \frac{\pi}{N}} <
\sqrt{\frac{N(N-1)}{2}}, \qquad \hbox{for all $N \geq 4$}.
\end{equation*}
We set $x=\pi /N$, so that the last inequality reads
\begin{equation*}
\frac{\sin x}{1-\cos x} < \sqrt{\frac{\pi(\pi-x)}{2x^2}}=\frac{1}{x}
\sqrt{\frac{\pi(\pi-x)}{2}}
\end{equation*}
We will now show that
\begin{equation*}
\frac{x \sin x}{1-\cos x} < \sqrt{\frac{\pi(\pi-x)}{2}} \qquad
\hbox{for all $x\in (0,\pi/4)$}.
\end{equation*}
Set
\begin{equation*}
f(x) = \frac{x\sin x}{1-\cos x}, \qquad g(x) =
\sqrt{\frac{\pi(\pi-x)}{2}}.
\end{equation*}
We now prove that $g-f$ is strictly decreasing in the interval
$(0,\pi/4)$ and that $g(\pi/4)-f(\pi/4)>0$. More precisely, we claim
that
\begin{equation*}
g'(x)-f'(x)= \frac{\operatorname{cosec}^2
\frac{x}{2}}{8\sqrt{\pi-x}} \left\{ \sqrt{2\pi} (\cos x-1) + 4
(x-\sin x) \sqrt{\pi-x} \right\}<0
\end{equation*}
All we have to show is that
\begin{equation*}
\sqrt{2\pi} (\cos x-1) + 4 (x-\sin x) \sqrt{\pi-x} < 0
\end{equation*}
for all $x\in (0,\pi/4)$, or
\begin{equation*} 
\sqrt{2\pi}(\cos x-1) < 4 \sqrt{\pi-x} (\sin x-x).
\end{equation*}
But
\begin{equation*}
\sin x -x > -\frac{x^3}{6}, \quad \cos x-1 <
-\frac{x^2}{2}+\frac{x^4}{4!}
\end{equation*}
and this implies that we can check the inequality
\begin{equation*}
\sqrt{2\pi} x^2 \left( \frac{x^2}{4!} - \frac12 \right) < -\frac{2}{3}
x^3 \sqrt{\pi -x},
\end{equation*}
or
\begin{equation*}
\sqrt{2\pi} (x^2-12) < -16x\sqrt{\pi-x}.
\end{equation*}
Both sides of this inequality are negative, since $x\in (0,\pi/4)$.
We now square and reverse the sense of the inequality, and have to
prove that
\begin{equation*}
2 \pi x^4 + 256 x^3 - 304\pi x^2 + 288\pi>0.
\end{equation*}
This is clearly true, since $2 \pi x^4 + 256 x^3 - 304\pi x^2 + 288\pi
> 2 \pi x^4 - 304\pi x^2 + 288\pi$ and this bi-quadratic equation has
no real roots in $(0,\pi/4)$. To complete the proof, we need to show
that $g(\pi/4)-f(\pi/4)>0$. But
\begin{equation*}
g(\pi/4)-f(\pi/4) = \frac{\pi}{2} \sqrt{\frac{3}{2}} -
\frac{\pi}{4\sqrt{2}-1}>0.
\end{equation*}
\end{proof}
\begin{lemma} \label{lem:1}
  For each fixed $N \geq 4$, the map $\Psi_N \colon [0,2] \to \R$ is
  strictly increasing.
\end{lemma}
\begin{proof}
With the same arguments contained in the proof of Lemma
\ref{lem:scemo}, one can prove that
\begin{equation*}
\frac{2
\tilde{r}_{12}^{-\al-2}-\tilde{r}_{13}^{-\al-2}}{\sum_{k=2}^4
\tilde{r}_{1k}^{-\al}} \quad \hbox{and} \quad
\frac{\tilde{r}_{12}^{-\al-2}}{\sum_{k=2}^N
\tilde{r}_{1k}^{-\al}} \qquad (N \geq 5)
\end{equation*}
are monotone functions of $\al \in (0,2)$. In particular for the
first one, we exploit the fact that $\tilde{r}_{12}=\tilde{r}_{14}$.
\end{proof}
We conclude using Lemma \ref{lem:phi}.
\begin{lemma} \label{lem:2}
For each $N \geq 4$, there results $\Psi_N(0) > 9/8$.
\end{lemma}
\begin{proof}
When $N=4$ we compute easily that $\Psi_4 (0) > 5/4 > 9/8$.  When
$N=5$ we use $\sin (\pi /5) = \sqrt{10-2\sqrt{5}}/4$ and easily
verify that $\Psi_5(0) > \frac{6\sqrt{5}-1}{5(\sqrt{5}-1)} > 9/8$.
For $N \geq 6$, we have to verify that $\Psi_N(0) = \Phi_N(0) + \frac{1}{2(N-1)} \frac{1}{\sin^2(\pi/N)}$.
From \eqref{eq:28} we have that this inequality is implied by
\begin{equation*}
\frac{1}{\sin^2 (\frac{\pi}{N})} > \frac{(N+8)(N-1)}{4N} \qquad
\hbox{for all $N \geq 6$}.
\end{equation*}
To conclude, recall that $\sin \frac{\pi}{N} < \frac{\pi}{N}$. Hence
\begin{equation*}
\frac{1}{\sin^2 (\frac{\pi}{N})} > \frac{N^2}{\pi^2},
\end{equation*}
and it is simple to check that for all $N \geq 6$ there results
\begin{equation*}
4N^3 > \pi^2 N^2 + 7\pi^2 N - 8\pi^2.
\end{equation*}
\end{proof}
\begin{theorem}
For every $N \geq 4$, there exists $\al_N < 1$ such that
\eqref{eq:36} holds for every $\al \geq \al_N$.
\end{theorem}
\begin{proof}
Let $f(\al) = \frac{(\al+2)^2}{8\al}$ be the right-hand
side of \eqref{eq:36}. This is a strictly decreasing function on the
interval $(0,2)$. Since $f(1)=9/8$, we deduce from lemmas
\ref{lem:1} and \ref{lem:2} that the graphs of $\Phi_N$ and of $f$
must intersect at a unique point $\al_N<1$. This concludes the proof.
\end{proof}


When $N$ is even and greater than $4$, we provide an example of a vertical variation that
verifies \eqref{eq:25} and also satisfies the \textit{hip-hop} symmetry (in $\RR^3$). 
Denoting $x_k=(\xi_k,z_k)\in\RR^3=\CC\times\RR$ the position of the $k$-th body,
this symmetry constraint, studied in \cite{cv} and \cite{tv}, imposes that
\[
 \forall k=1,\ldots,N-1 \qquad \xi_{k+1}(t)=e^{\frac{\pi i}{N}}\xi_k(t), \quad z_{k+1}(t)=-z_{k}(t).
\]
The simple variation on the $N$-gon configuration
introduced in \eqref{eq:w} is then no longer admissible as soon as $N \geq 6$; however, we can
consider the equivariant vector $w$ orthogonal to the plane of the
central configuration whose ``vertical'' components are $w_i=(-1)^i
/\sqrt{N}$, for~$i=1,2,\dots,N$.  Inequality~\eqref{eq:25bis} is then
equivalent to
\begin{equation}
\label{eq:comp1}
\frac{1}{N}\sum_{i,j}(-1)^{i+j}a_{ij} > h_N(\al), \qquad
\hbox{where} \quad h_N(\al)= \frac{(\al+2)^2}{8\al}
\frac{N}{2} \sum_{k=2}^N \frac{1}{r_{1k}^\al}.
\end{equation}
We already know (see Remark~\ref{rem:dueindici}) that whenever we
choose two consecutive bodies of the polygon (say the $i$-th and the
$(i+1)$-th) and we take $w_i=-w_{i+1}=1/\sqrt{2}$ then
\begin{equation*}
\frac{1}{2}\left(a_{ii} + a_{i+1\,i+1} -2a_{i\,i+1}\right) > h_N(\al).
\end{equation*}
In particular if we take $i=1,3,5, \ldots , N-1$ and we sum the
corresponding $N/2$ inequalities we obtain
\begin{equation*}
\frac{1}{2}\left(\sum_{i=1}^N a_{ii} -2\sum_{i \in
\{1,3,5,\dots,N-1\}} a_{i\,i+1} \right) > \tfrac{N}{2} h_N(\al),
\end{equation*}
or
\begin{equation}
\label{eq:comp2}
\frac{1}{N}\left(\sum_{i=1}^N a_{ii}
-2\sum_{i=1,3,5,\ldots,N-1}a_{i\,i+1} \right) = \frac{1}{N}\left( N
a_{11} -2 \frac{N}{2} a_{1\,2} \right) > h_N(\al).
\end{equation}
Comparing~\eqref{eq:comp1} and~\eqref{eq:comp2} we conclude if we can
prove
\begin{equation*}
-2\frac{N}{2}a_{1\,2} - 2N \sum_{j=3}^{N/2} (-1)^j a_{1\,j} -
 2\frac{N}{2} (-1)^{N/2 + 1} a_{1\,N/2 + 1} > 0
\end{equation*}
that is
\begin{equation*}
-a_{1\,2} - 2 \sum_{j=3}^{N/2} (-1)^j a_{1\,j} - (-1)^{N/2 + 1}
 a_{1\,N/2 + 1} > 0.
\end{equation*}
Replacing the expression of $a_{i\, j}$ (always negative when $i \neq
j$) we have to prove that
\begin{equation*}
\frac{1}{r_{1\,2}^{\al + 2}} + 2 \sum_{j=3}^{N/2}
                                \frac{(-1)^j}{r_{1\,j}^{\al + 2}} +
                                \frac{(-1)^{N/2 + 1}}{r_{1\,N/2 +
                                1}^{\al + 2}} > 0
\end{equation*}
and recalling that $r_{ij}(N) = \frac{2}{\sqrt{N}} \sin \left(
\frac{|i-j|}{N}\pi\right)$, inequality \eqref{eq:comp1} will follow
from
\begin{equation}
\label{eq:comp5}
g(N,\al)=\frac{1}{\sin^{\al + 2}\left( \frac{\pi}{N}\right)} + 2
\sum_{j=3}^{N/2} \frac{(-1)^j}{\sin^{\al + 2}\left(
\frac{j-1}{N}\pi\right)} + (-1)^{N/2 + 1} > 0.
\end{equation}
where $g$ is defined on the product $\{6,8,10,\ldots\} \times [0,2]$.
We then would like to prove~\eqref{eq:comp5} for every $N\geq6$, $N$
even, in an interval of values of $\al$ containing $\al=1$.  The
sum of the first two terms of the function $g$ is always positive,
indeed
\begin{equation*}
\frac{1}{\sin^{\al + 2}\left( \frac{\pi}{N}\right)} -
\frac{2}{\sin^{\al + 2}\left( \frac{2\pi}{N}\right)} =
\frac{\left[2\cos\left( \frac{\pi}{N}\right)\right]^{\al +
    2}-2}{\sin^{\al + 2}\left( \frac{2\pi}{N}\right)} =
\left(\frac{2}{\sin\left( \frac{2\pi}{N}\right)}\right)^{\al + 2}
\left[ \cos^{\al + 2}\left(\frac{\pi}{N}\right)-\frac{1}{2^{\al
      + 1}}\right]
\end{equation*}
and for every $N$ greater then $6$
\begin{equation*}
\cos\left(\frac{\pi}{N}\right) > \cos\left(\frac{\pi}{4}\right) =
\frac{\sqrt{2}}{2} \geq 2^{-\frac{\al + 1}{\al+2}}, \qquad
\forall \al \in [0,2].
\end{equation*}
The remaining terms can be collected in pairs of the kind
\begin{equation*}
\frac{1}{\sin^{\al +2}\left( \frac{(j-1)\pi}{N} \right)} -
\frac{1}{\sin^{\al +2}\left(\frac{j\pi}{N}\right)}, \qquad
\mbox{with $j$ even, $4 \leq j \leq \frac{N}{2}$}
\end{equation*}
whose sum is positive being $\sin x$ an increasing function when $x
\in [0,\pi/2]$ and $\al$ positive.  Concerning the other terms in
$g$, two different situations can occur: when $N/2$ is even or when
$N/2$ is odd.  In the first case the last two terms of the function
$g$ are $2\sin^{-\al - 2}\left( \frac{(N-2)\pi}{2N} \right)$ and
$-1$ whose sum is strictly positive.  When $N/2$ is odd there is just
a positive remaining term in $g$ which is $+1$.

\section{Asymptotic minimality for the weak--force case}

Equation \eqref{eq:2555} in Section \ref{sec:appl} suggests that
there should exist minimal colliding motions for small values
of~$\al$.  For the reader's sake, we will use a
somehow more transparent notation to stress the dependence on the parameter
$\al$ by writing $U \sb \al$ and
$\action_\al$ instead of $U$ and $\action$ defined in \eqref{eq:pot_hom} and 
\eqref{eq:19} respectively. Similarly $h_\al = \frac12 |\dot x_\al|^2 - U_\al (x_\al)$
denotes the energy of $x_\al$, $\cal{C}_\al \subset \elli$ the set of central
configurations of $U_\al$ and we refer to \eqref{eq:dyn_syst}$_\al$ to recall
the dynamical system \eqref{eq:dyn_syst} with the $\al$-homogeneous potential
$U_\al$.

Throughout this section, we will consider \emph{total} one-collision
solutions (see Definition \ref{def:1cs}) $x \sb \al$ to
\eqref{eq:dyn_syst}$_\al$ with the following ``initial conditions"
independent of $\al$:
\begin{equation*} 
|\dot x_\al(0)| = |x_\al (0)| =1.
\end{equation*}
In particular, the function $W$ defined in \eqref{eq:W} is identically
zero. In Section \ref{sec:vs} we have seen that there exists a diffeomorphism
$\Psi$ that rewrites any non-trivial $x \in H^1([0,1],\R^{Nd})$ 
in the coordinates $(\rho,s)$ and changes the time 
$t \in [0,1]$ into a new time $\tau \in [0,+\infty)$. 
On the space of these new variables we will use the norm
\begin{equation*}
\| (\zeta,v)  \|_{\mathcal{H}}^2 = \int_0^{+\infty} \left( |\zeta' (\tau)|^2 +
\zeta(\tau)^2  \right) \, d\tau
+ \int_0^{+\infty} \left( |v' (\tau)|^2 + |v(\tau)|^2  \right) \, d\tau
\end{equation*}
With an abuse of notation, we will continue to write $\action_\al$ instead of
$\action_\al \circ \Psi$.
Since the function $W$ is identically
zero, the Euler-Lagrange equations \eqref{eq:2} and \eqref{eq:4} reduce to
\begin{equation} \label{eq:2b}
-\left( \frac{4}{2-\al} \right)^2 \rho'' + \rho \left[|s'|^2 +
2U_\al (s) \right]= -\beta h_\al\rho^{\beta-1},
\end{equation}
\begin{equation} \label{eq:4b} 
-2\rho \rho' s' - \rho^ 2 s'' + \rho^ 2 \nabla {U_\al}_{|\elli} (s) =
\rho^2|s'|^2 s.
\end{equation}
It is convenient to introduce some terminology. 
\begin{definition}
Let $x \in H^1([0,1],\R^{Nd})$ be given, and let $s$ be the second component of $\Psi(x)$.
A function $\varphi$ is \emph{compactly supported variation on} $(a,b) \subset (0,1)$ corresponding to $x$ if $\varphi \in T_x H^1([0,1],\R^{Nd}) \simeq H^1([0,1],\R^{Nd})$ and if $v \colon [0,1] \to T_s \mathcal{E}$, where $\Psi(\varphi)=(\zeta,v)$.
\end{definition}
Furthermore, if $\Delta\subset\RR^{Nd}$ is the set of collision configurations defined in \eqref{collision_set}, fixed $\delta > 0$ we term $\Delta^\delta$ its open $\delta$-neighborhood. 
The following lemma is crucial for the proofs of Theorems
\ref{thm:first-exampl-homogr}, \ref{thm:non-homogr-moti} and \ref{th:ultimo}.

\begin{lemma}\label{lem:gra2:unif}
Let $(x_\al)_{\al\in A}$, $A\subset (0,2)$, be a family of 
one-collision solutions for \eqref{eq:dyn_syst}$_\al$, $\al \in A$,
and let $(K_\al)_{\al\in A}$, $K_\al \subset {\cal C}_\al$, 
be the sets of central configurations to which $x_\al$ is asymptotic.
for every $\al \in A$, in the sense of Definition \ref{def:as}.
If
\begin{equation}\label{no_coll_uniforme}
K_\al \subset {\cal C}_\al \setminus \Delta^\delta, \qquad \forall \al \in A
\end{equation}
for some $\delta>0$, then there exists $m>0$, which does not depend on $\al$, 
such that
\begin{equation*}
\nabla^2 {U_\al}_{|\mathcal{E}} (s_\al)(v,v) \geq -\al m |v|^2, 
\qquad \forall v \in T_{s_\al}\elli
\end{equation*}
where $s_\al = x_\al/|x_\al|$, for every $\al \in A$.
\end{lemma}

\begin{proof}
The conclusion follows from equations \eqref{eq:nabla2_U_centr} and
\eqref{eq:pot_hom_secondo} and the uniform assumption \eqref{no_coll_uniforme}.
\end{proof}

We recall some terminology.
\begin{definition}
A solution $x_\al(t)=(x_{\al,1},\ldots,x_{\al,N})(t)$ of
\eqref{eq:dyn_syst}$_\al$ is called \emph{homographic} 
if there exist two functions $\lambda(t)>0$ and $A(t)\in SO(d)$ 
and a fixed configuration $\bar s_\alpha = (\bar s_{\al,1},\ldots,\bar s_{\al,N})$ such that
\[
 x_\al(t)=\lambda(t)(A(t)\bar s_{\al,1},\ldots,A(t) \bar s_{\al,N})(t).
\]
\end{definition}
It is well known that the configuration $s_\al$ associated to a homographic motion is central, 
that is $s_\al$ is a critical point for the potential $U$ constrained to the ellipsoid $\elli$. See Definition \ref{def:cc} and Remark \ref{rem:cc}.

Let us state and prove the first result of this section.
\begin{theorem} \label{thm:first-exampl-homogr}
Let $(x_\al)_{\al\in A}$, $A\subset (0,2)$,
be a family of homographic one-collision solutions for \eqref{eq:dyn_syst}$_\al$
and let $\bar s_\al$ be such that
$x_\al(t) = r_\al(t) \bar s_\al$, for all $t \in [0,1]$ and $\al \in A$.
If there exists $\delta>0$ such that
\begin{equation} \label{eq:52}
(\bar s_\al)_\al \subset {\cal C}_\al \setminus \Delta^\delta
\end{equation}
then there exists $\bar \al \in (0,2)$ such that for every 
$\al \in (0,\bar\al)\cap A$ there exists $t_\al=t(x_\al)$ such that\footnote{With a slight abuse of notation, we write here $x_\alpha+\varphi$. This is justified by the fact that $\varphi$ is a tangent vector at $x_\alpha$ to the linear space $H^1([0,1],\R^{Nd})$.}
\begin{equation} \label{eq:delta}
\Delta \action_\al (x_\al,\varphi) := \action_\al (x_\al + \varphi) -
\action_\al (x_\al) > 0
\end{equation}
for every  compactly supported variation $\varphi$ on $(t_\al,1)$.

Furthermore if for every $\al\in A$ the energy of the homographic motion is
positive then there exists $\bar{\bar\al} \in (0,2)$ such that for every $\al
\in (0,\bar{\bar\al})\cap A$ inequality \eqref{eq:delta} holds
for every compactly supported variation $\varphi$ on $(0,1)$.
\end{theorem}
\begin{proof}
We will prove \eqref{eq:delta} by switching to the new coordinates $(\rho,s)$ and to the scaled time $\tau$ defined in
\eqref{time_riscal}. Since $\varphi$ is a compactly supported variation and by virtue of \eqref{eq:52}, we see that $\action_\al =\action_\al \circ \Psi$ is
smooth enough to write the Taylor expansion (where $(\zeta,v) = \Psi(\varphi)$
and $\supp v \subset (0,+\infty)$)
\begin{equation*}
\Delta \action_\al (x_\al,\varphi) = d^2 \action_\al (x_\al) \left(
(\zeta,v),(\zeta,v) \right) + o(\| (\zeta,v)  \|_{\mathcal{H}}^2),
\end{equation*}
where the first order term disappears because $x_\al$ is a critical point of $\action_\al$. 
Therefore, it is sufficient to prove that 
$d^2 \action_\al (x_\al)\left( (\zeta,v),(\zeta,v) \right) > 0$
whenever $\al$ is small enough and the support of $v$  is sufficiently away from
0. Equivalently, we will prove that 
\begin{equation} \label{eq:1} 
\frac{\partial^2 \action_\al}{\partial \rho^2} (\zeta,\zeta) + 2 
\frac{\partial^2 \action_\al}{\partial \rho \partial s} (\zeta,v)
+ \frac{\partial^2 \action_\al}{\partial s^2} (v,v) > 0,
\end{equation}
where 
\begin{eqnarray}
\frac{\partial^2 \action_\al}{\partial \rho^2} (\rho_\al,\bar
s_\al)(\zeta,\zeta) &=& \int_0^{+\infty} \left( \frac{4}{2-\al} \right)^2
(\zeta')^2 + \zeta^2 \left( |\bar s_\al'|^2 + 2 {U_\al}(\bar s_\al) \right) \, d
\tau, \label{eq:59} \\
\frac{\partial^2 \action_\al}{\partial \rho \partial s}(\rho_\al,\bar s_\al)
(\zeta,v) &=& 2 \int_0^{+\infty} \rho_\al \zeta \left[ \bar s_\al' \cdot v' +
\nabla {U_\al}_{|\mathcal{E}} (\bar s_\al) \cdot v \right] \, d\tau,
\label{eq:60} \\
\frac{\partial^2 \action_\al}{\partial s^2} (\rho_\al,\bar s_\al) (v,v) &=&
\int_0^{+\infty} \rho_\al^2 \left[ |v'|^2 + \nabla^2 {U_\al}_{|\mathcal{E}}
(\bar s_\al)(v,v) \right] \, d\tau. \label{eq:61}
\end{eqnarray}
Since we are dealing with homographic motions, $\bar s_\al$ is a constant
function, critical point of $U_\al$ constrained to the ellipsoid $\elli$, then 
$\displaystyle\frac{\partial^2 \action_\al}{\partial \rho \partial s}
(\rho_\al,\bar s_\al) (\zeta,v)=0$ for every pair $(\zeta,v)$.

As in the proof of Theorem \ref{th:main}, we introduce the auxiliary variable $w
= \rho_\al v$ and we wish to prove that
\begin{equation*}
\begin{split}
&\frac{\partial^2\action_\al}{\partial \rho^2}(\rho_\al,\bar s_\al)
(\zeta,\zeta) +
\frac{\partial^2 \action_\al}{\partial s^2}(\rho_\al,\bar s_\al)(v,v)\\
&\quad=\int_0^{+\infty} \left\{\left( \frac{4}{2-\al} \right)^2 (\zeta')^2 + 2
{U_\al}(\bar s_\al) \zeta^2 + \rho_\al^2 \left[|v'|^2+\nabla^2
{U_\al}_{|\mathcal{E}}(\bar s_\al)(v,v) \right]\right\} \, d\tau\\
&\quad=\int_0^{+\infty} \left\{\left( \frac{4}{2-\al} \right)^2 (\zeta')^2 +  2
{U_\al}(\bar s_\al)\zeta^2 + |w'|^2 +\frac{\rho_\al''}{\rho_\al}|w|^2 +\nabla^2
{U_\al}_{|\mathcal{E}}(\bar s_\al)(w,w)\right\} \, d\tau >0,
\end{split}
\end{equation*}
where in the last step we have integrated by parts.
Using the Euler-Lagrange equation \eqref{eq:2b} divided by $\rho$
we have
\begin{equation*}
\begin{split}
&\frac{\partial^2\action_\al}{\partial \rho^2}(\rho_\al,\bar s_\al)
(\zeta,\zeta) +
\frac{\partial^2 \action_\al}{\partial s^2}(\rho_\al,\bar s_\al)(v,v)\\
&\quad=\int_0^{+\infty} \left\{\left(\frac{4}{2-\al} \right)^2 (\zeta')^2 +
|w'|^2 
+ 2{U_\al} (\bar s_\al)\left[\zeta^2 + \left(\frac{2-\al}{4} \right)^2
|w|^2\right]\right\}\, d\tau \\ 
& \qquad\qquad+ \int_0^{+\infty} \left\{\left(\frac{2-\al}{4} \right)^2 \beta
h_\al \rho_\al^{\beta -2}|w|^2 + \nabla^2 {U_\al}_{|\mathcal{E}} (\bar s_\al)
(w,w) \right\} d\tau.
\end{split}
\end{equation*}
Since $U_\alpha$ is positive, by the uniform assumption \eqref{eq:52} and Lemma \ref{lem:gra2:unif}
there exist two positive constants $C,m$ such that
\begin{multline*}
\frac{\partial^2 \action_\al}{\partial \rho^2}(\rho_\al,\bar s_\al)
(\zeta,\zeta) +
\frac{\partial^2 \action_\al}{\partial s^2}(\rho_\al,\bar s_\al) (v,v)
\geq 
C\int_0^{+\infty} \left[ (\zeta')^2 +\zeta^2 + |w'|^2 \right] \,d\tau 
- \al m \int_0^{+\infty} |w|^2 \,d\tau\\
\quad + \left(\frac{2-\al}{4} \right)^2\int_0^{+\infty}  \beta h_\al
\rho_\al^{\beta -2}|w|^2\,d\tau.
\end{multline*}
If we suppose $h_\al>0$, for every $\al \in A$, we have that there exists
$C_1>0$ such that, whenever $\al$ is sufficiently small
\begin{equation*}
\begin{split}
\frac{\partial^2 \action_\al}{\partial \rho^2}(\rho_\al,\bar s_\al)
(\zeta,\zeta) +
\frac{\partial^2 \action_\al}{\partial s^2}(\rho_\al,\bar s_\al)(v,v)
& \geq \int_0^{+\infty} \left[ (\zeta')^2 +\zeta^2 + |w'|^2 \right] \,d\tau 
- \al M \int_0^{+\infty} |w|^2 \,d\tau\\
& \geq C_1 \| (\zeta,w) \|_{\mathcal{H}}^2
\end{split}
\end{equation*}
independently on the support of the function $w$.  Otherwise, when we do not
impose any assumption on the energy $h_\al$, since, for every $\al$, the
function $\rho_\al$ tends to 0 decreasing, we can find $\tau_\al$ sufficiently
large, such that
\begin{equation}\label{eq:123}
  \frac{\partial^2 \action_\al}{\partial \rho^2}(\rho_\al,\bar
s_\al)(\zeta,\zeta) +
  \frac{\partial^2 \action_\al}{\partial s^2}(\rho_\al,\bar s_\al)(v,v)
  \geq C_2 \| (\zeta,w) \|_{\mathcal{H}}^2
\end{equation}
for some positive constant $C_2$, whenever the support of $w$ is
contained in $(\tau_\al,+\infty)$.
\end{proof}
The asymptotic behavior of a collision solution, recalled in Proposition
\ref{prop:asymptotic}, suggests an extention of Theorem
\ref{thm:first-exampl-homogr} to suitable families of collision motions.
Also  in this case a uniform condition on the asymptotic sets of central
configuration will be assumed. With this aim we give the following definition.

\begin{definition}
  We say that the set of central configurations $K_\al \subset {\cal C}_\al$
  \emph{has the property of the asymptotic minimality} if for every $x_\al$
  solution \eqref{eq:dyn_syst}$_\al$ asymptotic to the set $K_\al$there
  exists $t_\al=t(x_\al)$ such that \eqref{eq:delta} holds for every $\varphi$
  with compact support in $(t_\al,1)$.
\end{definition}
The next result is a generalization of Theorem \ref{thm:first-exampl-homogr} to a larger class of total collision motions.
\begin{theorem} \label{thm:non-homogr-moti} Let $(K_\al)_{\al\in A}$, $A
  \subset (0,2)$, be a family of sets of central configurations. If there
  exists $\delta>0$ such that \eqref{no_coll_uniforme} holds then there exists
  $\bar \al \in (0,2)$ such that $K_\al$ has the property of the asymptotic
  minimality for every $\al \in (0,\bar \al)\cap A$.
\end{theorem}
\begin{proof}
  Let $x_\al$, for some $\al \in A$ be a solution of \eqref{eq:dyn_syst}$_\al$
  asymptotic to the set $K_\al$.We argue as in the proof of Theorem
  \ref{thm:first-exampl-homogr} and in this setting, replacing $w=\rho v$ and
using \eqref{eq:rhovi}, we integrate by parts and obtain
\begin{align*}
  &\frac{\partial^2 \action_\al}{\partial \rho^2}(\rho_\al,s_\al)(\zeta,\zeta)
  +
  \frac{\partial^2 \action_\al}{\partial s^2}(\rho_\al,s_\al)(v,v)\\
  &\quad=\int_0^{+\infty} \left\{\left( \frac{4}{2-\al} \right)^2 (\zeta')^2 +
    \zeta^2\left[|s_\al'|^2 + 2U_\al(s_\al)\right] + \rho_\al^2
\left[|v'|^2+\nabla^2 {U_\al}_{|\mathcal{E}}(s_\al)(v,v) \right]\right\} \,
d\tau\\
  &\quad=\int_0^{+\infty} \left\{\left( \frac{4}{2-\al} \right)^2 (\zeta')^2 +
    \zeta^2\left[|s_\al'|^2 + 2U_\al(s_\al)\right] + |w'|^2
+\frac{\rho_\al''}{\rho_\al}|w|^2 +\nabla^2
{U_\al}_{|\mathcal{E}}(s_\al)(w,w)\right\} \, d\tau\\
  &\quad=\int_0^{+\infty} \left\{\left(\frac{4}{2-\al} \right)^2 (\zeta')^2 +
    |w'|^2 + \left[\zeta^2 + \left( \frac{2-\al}{4} \right)^2 |w|^2\right]
    \left[|s_\al'|^2 + 2U_\al(s_\al)\right]\right\}\, d\tau \\
  &\qquad\qquad+ \int_0^{+\infty} \left\{\left(\frac{2-\al}{4} \right)^2 \beta
    h_\al \rho_\al^{\beta -2}|w|^2 + \nabla^2
    {U_\al}_{|\mathcal{E}}(s_\al)(w,w) \right\} d\tau .
\end{align*}
Since $|s'_\alpha|^2+2U_\alpha (s_\alpha)>0$, the uniform assumption \eqref{no_coll_uniforme} and Lemma
\ref{lem:gra2:unif} still imply the existence of a positive constant $C_2$
such that inequality \eqref{eq:123} holds.  On the other hand, integrating by parts the mixed
term $\frac{\partial^2 \action_\al}{\partial \rho \partial s} (\rho,s)
(\zeta,v)$ and recalling that $v$ has compact support we obtain
\begin{eqnarray}
\frac12 \frac{\partial^2 \action_\al}{\partial\rho\partial s}(\rho_\al,s_\al)
(\zeta,v) &=& \int_0^{+\infty} \rho_\al\zeta \left(s_\al'\cdot v' + \nabla
{U_\al}_{|\elli}(s_\al)\cdot v \right) \, d\tau \notag \\
&=& \int_0^{+\infty} \zeta s_\al' \cdot \left( w' -
\frac{\rho_\al'}{\rho_\al}w \right) + \zeta  \nabla {U_\al}_{|\elli}(s_\al) \cdot
w \, d\tau \notag \\
&=& \int_0^{+\infty} \left( -\zeta' s_\al' - \zeta s_\al'' - \zeta
\frac{\rho_\al'}{\rho_\al}s_\al' + \zeta \nabla {U_\al}_{|\elli}(s_\al)\right)
\cdot w \, d\tau. \label{eq:71bis}
\end{eqnarray}
Replacing the second Euler--Lagrange equation~\eqref{eq:4b} divided by $\rho^2$,
that is 
\begin{equation*}
-2 \frac{\rho'_\al}{\rho_\al} s'_\al -  s''_\al +\nabla {U_\al}_{|\elli}(s_\al) = |s'_\al|^2 s_\al,
\end{equation*}
into \eqref{eq:71bis} we obtain
\begin{equation*} 
\frac12\frac{\partial^2 \action_\al}{\partial\rho\partial s}(\rho_\al,s_\al)
(\zeta,v) 
= \int_0^{+\infty} \left( -\zeta' s_\al'+ \zeta \frac{\rho_\al'}{\rho_\al}s_\al'
+ \zeta |s_\al'|^2 s_\al\right) \cdot w \, d \tau.
\end{equation*}
The H\"{o}lder inequality gives immediately
\begin{eqnarray*}
\left| \int_0^{+\infty} \zeta' s_\al' \cdot w \, d \tau\right|
&\leq& \left(\int_0^{+\infty} |s_\al'|^2\, d \tau \right)^{1/2}
\left( \int_0^{+\infty} (\zeta')^2 |w|^2 \, d\tau \right)^{1/2} \\
&\leq& \|w\|_\infty \left(\int_0^{+\infty} |s_\al'|^2\, d \tau \right)^{1/2}
\left( \int_0^{+\infty} (\zeta')^2\, d\tau \right)^{1/2}
\\
\left| \int_0^{+\infty} \zeta \frac{\rho'_\alpha}{\rho_\alpha} s_\al' \cdot w \, d \tau\right|
&\leq& \left(\int_0^{+\infty} |s_\al'|^2\, d \tau \right)^{1/2}
\left( \int_0^{+\infty} (\zeta)^2 \left( \frac{\rho'_\alpha}{\rho_\alpha} \right)^2 |w|^2 \, d\tau \right)^{1/2} \\
&\leq& \|w\|_\infty \left(\int_0^{+\infty} |s_\al'|^2\, d \tau \right)^{1/2}
\left( \int_0^{+\infty} (\zeta)^2 \left( \frac{\rho'_\alpha}{\rho_\alpha} \right)^2 \, d\tau \right)^{1/2} \\
\left| \int_0^{+\infty} \zeta |s_\al'|^2 s_\al\cdot w \, d \tau\right| 
&\leq& \|\zeta\|_\infty
\|w\|_\infty \int_0^{+\infty} |s_\al'|^2 \, d \tau.
\end{eqnarray*}
and then
\begin{multline*}
\frac12\left|\frac{\partial^2 \action_\al}{\partial\rho\partial s}
(\rho_\al,s_\al)(\zeta,v)\right| \leq \\
 \|w\|_\infty \left[ \left(\int_0^{+\infty} (\zeta')^2\,d\tau\right)^{\frac12}
+ \left(\int_0^{+\infty} \zeta^2
\left(\frac{\rho'_\al}{\rho_\al} \right)^2\, d \tau  \right)^{\frac12} \right] \left(\int_0^{+\infty}|s_\al'|^2\,d\tau\right)^{\frac12}
+ \|\zeta\|_\infty\int_0^{+\infty}|s_\al'|^2\,d\tau.
\end{multline*}
Recall that $\al$ is fixed (and so small that \eqref{eq:123}
holds). Pick now $\varepsilon >0$. Since $\rho_\al' / \rho_\al$
converges to a finite limit as $\tau \to {+\infty}$, and
$\int_0^{+\infty} |s_\al'|^2 < {+\infty}$ (Proposition \ref{prop:asymptotic}),
there exists $\tau_\al$ depending on $x_\al$ such that  
\begin{equation*}
\int_0^{+\infty} |s_\al'|^2\, d \tau = \int_{\supp w} |s'_\al|^2\, d \tau <
\varepsilon
\end{equation*}
whenever $\supp w \subset (\tau_\al,{+\infty})$.
Hence, for all such $w$'s,
\begin{equation}\label{eq:122}
\left| \frac{\partial^2 \action_\al}{\partial \rho \partial s}(\rho_\al,s_\al)
(\zeta,v) \right| 
\leq C_3 \sqrt{\varepsilon} \ \|w\|_\infty \, \|\zeta\|_{H^1} 
\leq C_4 \sqrt{\varepsilon}\  \|(\zeta,w) \|_{\mathcal{H}}^2
\end{equation}
for some positive constants $C_3,C_4$.
From \eqref{eq:123} and \eqref{eq:122} we obtain
\begin{equation*}
d^2 \action_\al (x_\al) ((\zeta,v),(\zeta,v)) \geq \left( C_2 - 2C_4
\sqrt{\varepsilon} \right) \|(\zeta,w) \|_{\mathcal{H}}^2 > 0.
\end{equation*}
\end{proof}

In Theorem \ref{thm:non-homogr-moti}, we cannot exclude that, as
$\alpha \to 0$, the support of the variation $\varphi$ moves off to
the collision time $t=1$. It is natural to investigate under what
circumstances it is possible to single out a time $t^*$, independent
of $\alpha$, such that the second differential $d^2 \action_\alpha
(x_\alpha)(\varphi,\varphi)>0$ is positive for any variation $\varphi$
supported in $(t^*,1)$. It will turn out that the following uniform condition
on the behavior of $\rho_\alpha$ plays a crucial r\^ole.

\bigskip

\begin{description}
 \item[(UC)] As $\tau \to +\infty$, $\rho_\alpha(\tau) \to 0$ uniformly with
respect to
$\alpha \in (0,2)$. More precisely, for all $\sigma >0$ there exists
$\tau_\sigma >0$ such that for all $\alpha \in (0,1)$ and all $\tau \geq
\tau_\sigma$ there results $\rho_\alpha (\tau) < \sigma$.
\end{description}

\bigskip

Since it is clear that it would be useless to take the limit as
$\alpha \to 0$ inside (\ref{eq:dyn_syst}), we need to single out a
non-trivial \textit{limiting problem} that describes the asymptotic
properties of one-collision solutions.  Therefore we introduce the
scaled potential
\begin{equation} \label{eq:10}
\widetilde U_\al(x) = \frac{1}{\al} \sum_{i<j} 
\frac{m_i m_j }{|x_i-x_j|^{\al}}=
\frac{1}{\al}U_\al(x).
\end{equation}
The corresponding action reads
\begin{equation} \label{eq:18}
\widetilde{\action}_\alpha (x) = \frac{1}{2} \int_0^1 |\dot x|^2 \, dt +
\frac{1}{\alpha} \int_0^1 U_\alpha (x) \, dt 
\end{equation}
When we replace $\tilde U_\al$ to $U_\al$ in \eqref{eq:dyn_syst}$_\al$, 
the solutions of the new
dynamical system are strictly linked to the solutions of the old one as the next Lemma
asserts. Its very simple proof is omitted.
\begin{lemma}
If $\tilde{x} = \al^{-\frac{1}{\al+2}} x$, then
$\widetilde{\action}_\al (\tilde{x}) = \al^{-\frac{2}{\al+2}}
\action_\al (x)$. In particular, if $x_\al$ is a solution of
\eqref{eq:dyn_syst}$_\al$, then $\tilde{x}_\al =
\al^{-\frac{1}{\al+2}} x_\al$ solves
\begin{equation} \label{eq:13}
\ddot {\tilde{x}}_\alpha = \nabla \widetilde{U}_\al (\tilde{x}_\al).
\end{equation}
\end{lemma}

\begin{remark} \label{rem:6.8}
It is evident that a solution of $\ddot x = \nabla \widetilde{U}_\alpha(x)$ is
also
a solution of $\ddot x = \nabla \widehat{U}_\alpha(x)$, where
\begin{equation} \label{eq:11}
\widehat{U}_\al(x) = \frac{1}{\al} \sum_{i<j} m_i m_j \left(
\frac{1}{|x_i-x_j|^{\al}} -1 \right).
\end{equation}
For each $x=(x_1,\dots,x_n) \in \R^{Nd}$, there results
\begin{equation} \label{eq:12}
\lim_{\alpha \to 0} \widehat{U}_\al (x) = -\sum_{i<j} m_i m_j \log |x_i-x_j| =:
U_{\mathrm{log}} (x).
\end{equation}
However, the potential $\widehat{U}_\alpha$ is lacking the homogeneity property
which seems to be essential in the definition of the new variables $(\rho,s)$,
see Section~\ref{sec:vs}.
\end{remark}
We consider a family $(\tilde x_\al)_\al$, $\al\in(0,2)$, such that, fixed $\al$,
$\tilde{x}_\al$ solves
\begin{equation}
\begin{cases}
 \ddot {\tilde{x}}_\alpha = \nabla \widetilde{U}_\al (\tilde{x}_\al) \\
\tilde{x}_\alpha (0) = x_\alpha^0 \in \R^{Nd}, \\
\dot{\tilde{x}}_\alpha (0) = v_\alpha^0 \in \R^{Nd}.
\end{cases}
\end{equation}
We express $\tilde x_\al$ in terms of the generalized polar coordinates
$(\tilde{\rho}_\alpha,\tilde{s}_\alpha)$ and the new time $\tau$. Hence
$\tilde{\rho}_\alpha$ and $\tilde{s}_\alpha$ satisfy the Euler--Lagrange
equations \eqref{eq:2b} and \eqref{eq:4b} with $\widetilde{U}_\alpha$ instead of
$U_\alpha$. We make the following assumptions on the initial condition:
\begin{description}
\item[(IC1)] $\tilde{\rho}_\alpha (0)=1$ and $\tilde{\rho}'_\alpha (\tau)< 0$ 
for all $\alpha$ and for all $\tau \geq 0$.
 \item[(IC2)] $\tilde{s}_\alpha (0) \to s_0$ and $\tilde{s}'_\alpha (0) \to v_0$
as $\alpha \to 0$.
\end{description}
Define
\begin{equation} \label{eq:14}
\Gamma_\alpha (\tau) = \frac{1}{2} |\tilde s'_\alpha(\tau)|^2 - \widehat{U}_\alpha
(\tilde s_\alpha (\tau))
\end{equation}
where $\hat U_\al$ has been introduced in \eqref{eq:11}. 
In this setting we prove the next four Lemmas.
\begin{lemma} \label{lem:esplode0}
There exists a constant $C>0$, independent of $\alpha$, such
that
\begin{equation} \label{eq:15}
\int_0^{+\infty} -\frac{\tilde \rho'_\alpha}{\tilde \rho_\alpha} 
|\tilde s'_\alpha|^2 \, d \tau \leq C.
\end{equation}
\end{lemma}
\begin{proof}
By differentiating $\Gamma_\alpha$ (with respect to $\tau$) and making use of
\eqref{eq:4b} with $U_\al$ replaced by $\tilde U_\al$, we compute
\begin{equation}
\frac{d \Gamma_\alpha}{d \tau} = -2 \frac{\tilde \rho'_\alpha}{\tilde\rho_\alpha}
|\tilde s'_\alpha|^2.
\end{equation}
Therefore
\begin{equation} \label{eq:16}
\int_0^{+\infty} -\frac{\tilde\rho'_\alpha}{\tilde\rho_\alpha} |\tilde s'_\alpha|^2 \, d
\tau = \lim_{\tau \to +\infty} \Gamma_\alpha (\tilde s_\alpha (\tau)) - \Gamma_\alpha
(\tilde s_\alpha (0))
\end{equation}
We will complete the proof by showing that the right--hand side of \eqref{eq:16}
has a finite limit as $\alpha \to 0$. Since $\Gamma_\alpha (0) = \frac{1}{2}
|\tilde s'_\alpha(0)|^2 - \hat{U}_\alpha (\tilde s_\alpha (0))$, by virtue of
assumption \textbf{(IC2)}, we have that $\Gamma_\alpha (0)$ has a limit as $\alpha \to 0$.

As regards the behavior of $\Gamma_\alpha (+\infty) := \lim_{\tau \to +\infty}
\Gamma_\alpha (\tilde s_\alpha (\tau))$, we deduce from the asymptotic estimates (see
Proposition \ref{prop:asymptotic} (c)) that $\lim_{\tau \to +\infty}
\widehat{U}_\alpha (\tilde s_\alpha (\tau))=\hat{b}_\alpha$ exists and is finite. We
choose $\bar s_\alpha \in \mathcal{E}$ such that $\hat{b}_\alpha =
\widehat{U}_\alpha (\bar s_\alpha)$. Since $\mathcal{E}$ is a compact set, we
may assume that $\bar{s}_{\alpha_k} \to \bar{s}_0$ for a suitable subsequence
$\alpha_k \to 0$. Moreover, it is known that $\lim_{\tau \to +\infty}
|\tilde s'_{\alpha_k} (\tau)| =0$ (see Proposition \ref{prop:asymptotic} (f)).
We conclude as before that $\Gamma_{\alpha_k}
(+\infty)$ has a finite limit as $\alpha_k \to 0$.
\end{proof}

\begin{lemma}\label{lem:esplode1}
For any $\ge>0$ there exist $\tau_\ge >0$ and $\al_\ge \in(0,2)$ such that
\begin{equation}\label{eq:5}
\frac{2-\alpha}{4\alpha} \left( 1-
\tilde{\rho}_\alpha(\tau)^{4\alpha/(2-\alpha)} \right) \geq \frac{1}{\ge}
\end{equation}
for all $\alpha \in (0,\al_\ge)$ and $\tau \geq \tau_\ge$.
\end{lemma}
\begin{proof}
To save notation, we set $\gamma = \gamma(\alpha) = 4\alpha/(2-\alpha)$.
First of all, we remark that if $\rho_1$, $\rho_2\in (0,1)$ with
$\rho_1<\rho_2$, then
\begin{equation*}
\frac{1-\rho_1^\gamma}{\gamma} > \frac{1-\rho_2^\gamma}{\gamma}.
\end{equation*}
We fix $\ge >0$ and choose $\eta=\eta_\ge>0$ such that $-\log \eta >
1/\ge$. From assumption {\bf (UC)}, we can fix $\tau_\ge>0$ such
that $\tilde{\rho}_\alpha (\tau) < \eta$ whenever $\alpha \in (0,2)$ and $\tau
\geq \tau_\ge$. Furthermore, since as $\al \to 0$ $(1-\eta^\gamma)/\gamma \to -\log \eta$,
we can fix $\alpha_\ge \in (0,2)$ such that
\begin{equation*}
\frac{1-\eta\sp\gamma}{\gamma} \geq -\log \eta -\ge.
\end{equation*}
Finally, if $\tau \geq \tau_\ge$ and $\alpha < \alpha_\ge$, we get
\begin{equation} \label{eq:6}
\frac{1-\tilde{\rho}_\alpha(\tau)^\gamma}{\gamma} > \frac{1-\eta^\gamma}{\gamma}
\geq
-\log \eta -\ge \geq \frac{1}{\ge}-\ge.
\end{equation}
Inequality (\ref{eq:5}) is of course equivalent to (\ref{eq:6}).
\end{proof}
Before proceeding, we notice that, since $h_\alpha = \frac{1}{2}|\dot x_\al|^2 - U_\al
(x_\al)$ and $\tilde{x}_\al = \al^{-\frac{1}{\al+2}}{x}_\al$, there results 
$\tilde{h}_\al = \frac{1}{2} |\dot{\tilde{x}}_\al|^2-\frac{1}{\al}
\tilde{U}_\al(\tilde{x}_\al) = \al^{-\frac{1}{\al+2}} h_\al$. Similarly, from Remark
\ref{rem:6.8} we also get $\hat{h}_\al = \frac{1}{2} |\dot{\tilde{x}}_\al|^2
-\hat{U}_\al (\tilde{x}_\al) = \tilde{h}_\al + \frac{1}{\al} \sum_{i<j} m_i m_j$. 
We will assume that $\hat{h}_\al=0$, which amounts to

\medskip
\noindent\textbf{(H)} The energy of the solution $x_\alpha$ is $h_\alpha = \frac{\sum_{i<j}m_i m_j}{\alpha^{\frac{\alpha}{\alpha+2}}}$.
\medskip 

\begin{lemma} \label{lem:disotto}
Assume condition \textbf{(H)}. Then there exists a constant $C$ such that
\begin{equation*}
 2 \widetilde{U}_\alpha (\tilde{s}_\alpha (\tau)) + \beta \tilde h_\alpha
\tilde{\rho}_\alpha(\tau)^{\beta -2} \geq C + \frac{1}{\ge}
\end{equation*}
for all $\alpha \in (0,\al_\ge)$ and $\tau \geq \tau_\ge$.
\end{lemma}
\begin{proof}
We can write
\begin{multline*}
{\widetilde U_\al} (\tilde{s}_\al)+\frac{\beta}{2} \tilde h_\al \rho_\al^{\beta -2} =
\sum_{i<j}m_i m_j
\left[\frac{1}{\al}\frac{1}{|\tilde{s}_{\alpha,i}-\tilde{s}_{\alpha,j}|^\al}
-\frac{1}{\al}\frac{2+\al}{(2-\al)}
\tilde{\rho}_\al^{4\al/(2-\al)}\right]\\
= \sum_{i<j}m_i m_j
\left[\frac{1}{\al}\left(\frac{1}{|\tilde{s}_{\alpha,i}-\tilde{s}_{\alpha,j}
|^\al}-1\right)+\frac{1}{\al}
-\frac{1}{\al}\frac{2+\al}{(2-\al)}\left(\tilde{\rho}_\al^{4\al/(2-\al)}
-1\right)
-\frac{1}{\al}\frac{2+\al}{2-\al}\right]\\
= \sum_{i<j}m_i m_j
\left[\frac{1}{\al}\left(\frac{1}{|\tilde{s}_{\alpha,i}-\tilde{s}_{\alpha,j}
|^\al}-1\right)
+\frac{1}{\al}\left(1-\frac{2+\al}{2-\al}\right)
-4\frac{2+\al}{(2-\al)^2}\frac{2-\al}{4\al}\left(\tilde{\rho}_\al^{4\al/(2-\al)}
-1\right)\right].
\end{multline*}
We now observe that since $\tilde s_\al \in \elli$, then $|\tilde{s}_{\alpha,i}-\tilde{s}_{\alpha,j}|\leq 2$, for all $\al$ and $i \neq j$, and
\[
\frac{1}{\al}\left(\frac{1}{|\tilde{s}_{\alpha,i}-\tilde{s}_{\alpha,j}
|^\al}-1\right)\geq\frac{1-2^\al}{\al 2^\al}
\]
where the right hand side converges to $-\log 2$. 
The conclusion follows from Lemma \ref{lem:esplode1} and easy algebraic
inequalities.
\end{proof}
\begin{remark}
We notice that assumption \textbf{(H)} implies in particular $\lim_{\alpha \to 0} h_\alpha = \sum_{i<j}m_i m_j$. More generally, the same proof adapts to the case in which $\hat{h}_\alpha =C$, a constant independent of $\alpha$. Indeed, the ``old'' energy would be $h_\alpha = \left( C \alpha^{\frac{2}{\alpha+2}} - \frac{1}{\alpha^{\frac{\alpha}{\alpha+2}}} \right) \sum_{i<j} m_i m_j$, and the first term tends to zero as $\alpha \to 0$.
\end{remark}

\begin{lemma} \label{lem:esplode2}
Assume condition \textbf{(H)}. Then, for every $\ge >0$ there exists $\tau_\ge >0$ and $\al_\ge \in (0,2)$ such that
\begin{equation*}
\int_{\tau_\ge}^\infty |\tilde s'_\alpha|^2 \, d\tau < \ge.
\end{equation*}
for all $\alpha \in (0,\al_\ge)$.
\end{lemma}
\begin{proof}
Set $\phi_\alpha(\tau) = - \tilde{\rho}'_\alpha (\tau) / \tilde{\rho}_\alpha
(\tau)$. By direct computation 
\[
\phi'_\alpha(\tau) = - \frac{\tilde{\rho}''_\al}{\tilde{\rho}_\alpha} +
\left( \frac{\tilde{\rho}'_\alpha}{\tilde{\rho}'_\alpha} \right)^2 = 
-\frac{\tilde{\rho}''_\al}{\tilde{\rho}_\alpha} + \phi_\alpha (\tau)^2
\]
hence using equation \eqref{eq:2b} and Lemma \ref{lem:disotto} we have
\begin{eqnarray*}
\phi'_\alpha(\tau)
&=& \left( \frac{2-\alpha}{4} \right)^2 \left[ -| \tilde{s}'_\alpha|^2 - 2
\widetilde{U}_\alpha (\tilde{s}_\alpha) - \beta h_\alpha
\tilde{\rho}_\alpha^{\beta-2} \right] + \phi_\alpha (\tau)^2 \\
&\leq& - \left( C+\frac{1}{\ge} \right) + \phi_\alpha (\tau)^2.
\end{eqnarray*}
Since, by \textbf{(IC1)}, $\tilde{\rho}'_\alpha < 0$ we deduce that
$\phi_\alpha$ is positive. We claim that 
\begin{equation} \label{eq:17}
 \phi_\alpha(\tau)^2 \geq C+\frac{1}{\ge} \qquad\hbox{for every $\tau > 0$}.
\end{equation}
If \eqref{eq:17} is false, then $\phi_\alpha (\tau_0)< \sqrt{C+{1}/{\ge}}$ for
some $\tau_0>0$. Consider a solution $\psi$ of the Cauchy problem
\begin{equation*}
 \begin{cases}
  \psi' = \psi^2 -(C+\frac{1}{\ge}) \\
   \psi(\tau_0) = \phi_\alpha (\tau_0).
 \end{cases}
\end{equation*}
A basic comparison theorem for ODEs implies that $\phi_\alpha (\tau) \leq
\psi(\tau)$ for all $\tau \geq \tau_0$. But $\psi(\tau) \to -(C+1/\ge)$ as $\tau
\to +\infty$, and therefore $\phi_\alpha$ becomes negative for sufficiently
large times. This is a contradiction that proves \eqref{eq:17}.
\end{proof}
\begin{theorem} \label{th:ultimo}
Let $(x_\alpha)_{\alpha \in A}$, $A \subset (0,2)$, be a family of
total one-collision solutions of \eqref{eq:dyn_syst}$_\alpha$, and let
$(\tilde{x}_\alpha)_{\alpha \in A}$ be the corresponding solutions for the
potential $\widetilde{U}_\alpha$. Retain assumptions \textbf{(UC)},
\textbf{(IC1-2)} and \textbf{(H)}. If there exists $\delta >0$ for which
\eqref{no_coll_uniforme} holds, then there exist $t^*$ and a sequence 
$(\alpha_k)_k \subset A$ $\alpha_k \to 0$ such that
\begin{equation}
 d^2 \action_{\alpha_k} (x_{\alpha_k})(\varphi,\varphi) >0
\end{equation}
for every variation $\varphi$ with support in $(t^*,1)$.
\end{theorem}
\begin{proof}
As already remarked, we can consider the action $\widetilde{\action}_\alpha$
instead
of $\action_\alpha$. Furthermore, \eqref{eq:59}, \eqref{eq:60} and \eqref{eq:61}
hold with $U_\alpha$ replaced by $\widetilde{U}_\alpha$.
When we compute the variation of the action functional we follow
the proof of Theorem \ref{thm:non-homogr-moti} to obtain
\begin{equation*}
\begin{split}
&\frac{\partial^2 \widetilde\action_\al}{\partial \rho^2}
(\tilde \rho_\al,\tilde s_\al) (\zeta,\zeta) +
\frac{\partial^2 \widetilde\action_\al}{\partial \tilde{s}_\al^2}
(\tilde \rho_\al,\tilde s_\al)(v,v)\\
&\qquad \qquad =\int_0^{+\infty} \left\{\left(\frac{4}{2-\al} \right)^2 (\zeta')^2 + |w'|^2 
+ (\zeta)^2 \left[ |\tilde{s}_\al'|^2 + 2{\widetilde U_\al} (\tilde{s}_\al)
\right]\right\}\, d\tau \\ 
&\qquad + \int_0^{+\infty} \nabla^2 (\widetilde
U_\al)_{|\mathcal{E}}(\tilde{s}_\al)(w,w)\,d\tau
+ \int_0^{+\infty} 
\left(\frac{2-\al}{4} \right)^2 |w|^2\left[|\tilde{s}'_\al|^2 + 2{\widetilde
U_\al}
(\tilde{s}_\al)+\beta
h_\al \tilde{\rho}_\al^{\beta -2}  \right] d\tau.
\end{split}
\end{equation*}
The first integral is positive and tends to $+\infty$ as $\al \to 0$; the
second one is bounded from below (indeed, following the same idea of Lemma
\ref{lem:gra2:unif}, $\nabla^2 (\widetilde U_\al)_{|\mathcal{E}}(s)(w,w) \geq -m
|w|^2$).
The third one can be handled with the help of Lemma \ref{lem:disotto},
giving us the estimate
\[
\int_0^{+\infty} \left(\frac{2-\al}{4} \right)^2 |w|^2\left[|\tilde{s}'_\al|^2
+ 2{\widetilde U_\al} (\tilde{s}_\al)+\beta
\tilde h_\al \tilde{\rho}_\al^{\beta -2}  \right] d\tau \geq \left( K + \frac{1}{\ge}
\right) \int_0^{+\infty} |w|^2 \, d\tau,
\]
for some $K>0$.
In conclusion, there exists a constant $C_1>0$ such that 
\[
\frac{\partial^2 \tilde\action_\al}{\partial \rho^2}(\rho_\al,\bar
s_\al)(\zeta,\zeta) +
\frac{\partial^2 \tilde\action_\al}{\partial s^2}(\rho_\al,\bar s_\al)(v,v)
\geq C_2 \| (\zeta,w) \|_{\mathcal{H}}^2.
\]
Concerning the mixed derivative, we argue as in the proof of
\ref{thm:non-homogr-moti} to obtain
\begin{equation*} 
\frac12\frac{\partial^2 \widetilde{\action}_\al}{\partial \rho\partial
s}(\tilde{\rho}_\al,\tilde{s}_\al)
(\zeta,v) 
= \int_{\supp v} \left( -\zeta' \tilde{s}'_\al+ \zeta
\frac{\tilde{\rho}'_\al}{\tilde{\rho}_\al} \tilde{s}'_\al
+ \zeta | \tilde{s}'_\al|^2 \tilde{s}_\al\right) \cdot w \, d \tau.
\end{equation*}
We use the H\"{o}lder inequality and get
\begin{multline} \label{eq:74}
\frac12 \left| \frac{\partial^2 \widetilde{\action}_\al}{\partial \rho\partial
s}(\tilde{\rho}_\al,\tilde{s}_\al)(\zeta,v) \right|
\leq \left( \int_{\supp v} |\tilde{s}'_\alpha|^2  \, d\tau \right)^{\frac12}\times\\
\times\left[\|w\|_\infty\left( \int_{\supp v} (\zeta')^2  \, d\tau
\right)^{\frac12} 
+ \|v\|_\infty \|\zeta\|_\infty \left( \int_{\supp v} (\tilde{\rho}'_\al)^2 \,
d\tau \right)^{\frac12} \right. \\
\left.
{} + \|w\|_\infty\|\zeta\|_\infty\left( \int_{\supp v} |\tilde{s}'_\alpha|^2  \, d\tau \right)^{\frac12} \right] .
\end{multline}
We have seen in Lemma \ref{lem:esplode2} that, for every $\ge >0$, there exists
$\tau_\ge$, independent of $\alpha$, such that 
\[
 \int_{\supp v} |\tilde{s}'_\alpha|^2  \, d\tau  < \ge \qquad \hbox{for all
$\alpha$},
\]
provided that $\supp v \subset (\tau_\ge,+\infty)$.
To conclude, we notice that we are integrating over the compact set $\supp v$
which is disjoint from the collision time. Basic results in the theory of ODEs
(see Theorem 8.4 in \cite{amann})
imply that $x'_\alpha$ converges locally uniformly --- and locally in $L^2$ ---
to some limit as $\alpha \to 0$. Thus, also $\rho'_\alpha$ converges in
$L^2(\supp v)$ to a limit, and in particular
\[
 \sup_\alpha \int_{\supp v} (\tilde{\rho}'_\al)^2 \, d\tau < + \infty.
\]
This and \eqref{eq:74} give that
\[
 \frac12 \left| \frac{\partial^2 \widetilde{\action}_\al}{\partial \rho\partial
s}(\tilde{\rho}_\al,\tilde{s}_\al)(\zeta,v) \right| \leq C \sqrt{\ge} \
\|(\zeta,v)\|_{\mathcal{H}}^2.
\]

\end{proof}



%
\end{document}